\documentclass[aop,MSNbibl,seceqn,citesort,dvips]{arximspdf}

\doi{10.1214/11-AOP717} 
\volume{41}
\issue{3B}
\pubyear{2013}
\firstpage{2225}
\lastpage{2260}

\makeatletter

\newtheorem{proposition}{Proposition}[section]
\newtheorem{corollary}[proposition]{Corollary}
\newtheorem{theorem}[proposition]{Theorem}
\newtheorem{lemma}[proposition]{Lemma}

\newproclaim{remark}[proposition]{Remark}

\newcommand{\N}{\mathbf{N}}
\newcommand{\Z}{\mathbf{Z}}
\newcommand{\R}{\mathbf{R}}

\renewcommand{\P}{\mathrm{P}}
\newcommand{\E}{\mathrm{E}}
\newcommand{\F}{\mathcal{F}}

\newcommand{\one}{\mathbf{1}}
\renewcommand{\d}{{\mathrm d}}

\makeatother

\begin{document}
\begin{frontmatter}

\title{On the chaotic character of the stochastic heat equation,
before the onset of intermitttency}
\runtitle{On the chaotic character of the stochastic heat equation}

\begin{aug}
\author[A]{\fnms{Daniel} \snm{Conus}\corref{}\thanksref{t1}\ead[label=e1]{daniel.conus@lehigh.edu}},
\author[B]{\fnms{Mathew} \snm{Joseph}\thanksref{t2}\ead[label=e2]{joseph@math.utah.edu}}
\and
\author[B]{\fnms{Davar} \snm{Khoshnevisan}\thanksref{t3}\ead[label=e3]{davar@math.utah.edu}}
\runauthor{D. Conus, M. Joseph and D. Khoshnevisan}
\affiliation{Lehigh University, University of Utah and University of Utah}
\address[A]{D. Conus\\
Department of Mathematics\\
Lehigh University\\
Christmas--Saucon Hall\\
14 East Packer Avenue\\
Bethlehem, Pennsylvania 18015\\
USA\\
\printead{e1}} 
\address[B]{M. Joseph\\
D. Khoshnevisan\\
Department of Mathematics\\
University of Utah\\
155 South 1400 East\\
Salt Lake City, Utah 84112-0090\\
USA\\
\printead{e2}\\
\hphantom{E-mail: }\printead*{e3}}
\end{aug}

\thankstext{t1}{Supported in part by the
Swiss National Science Foundation Fellowship PBELP2-122879.}

\thankstext{t2}{Supported in part by the
NSF Grant DMS-07-47758.}

\thankstext{t3}{Supported in part by the
NSF Grants DMS-07-06728 and DMS-10-06903.}

\received{\smonth{3} \syear{2011}}
\revised{\smonth{11} \syear{2011}}

%
\begin{abstract}
We consider a nonlinear stochastic heat equation $\partial_t u =
\frac12 \partial_{xx} u + \sigma(u)\partial_{xt}W$, where
$\partial_{xt}W$ denotes space--time white noise and $\sigma\dvtx\R\to\R
$ is
Lipschitz continuous. We establish that, at every fixed time $t>0$, the
global behavior of the solution depends in a critical manner on the
structure of the initial function $u_0$: under suitable conditions on
$u_0$ and $\sigma$, $\sup_{x\in\R}u_t(x)$ is a.s. finite when $u_0$
has compact support, whereas with probability one,
$\limsup_{|x|\to\infty} u_t(x)/({\log}|x|)^{1/6}>0$ when $u_0$ is bounded
uniformly away from zero. This sensitivity to the initial data of the
stochastic heat equation is a way to state that the solution to the
stochastic heat equation is \textit{chaotic} at fixed times, well before
the onset of intermittency.
\end{abstract}

%

\end{frontmatter}

\section{Introduction and main results}\label{intro}

Let $\mathbf{W}:=\{W(t,x)\}_{t\ge0,x\in\R}$ denote a real-valued Brownian
sheet indexed by two parameters $(t,x)\in\R_+\times\R$.
That is, $\mathbf{W}$ is a centered Gaussian process with covariance
%
\begin{equation}
\operatorname{Cov}( W(t,x), W(s,y))
= \min(t,s)\times\min(|x|,|y|)\times\one_{(0,\infty)}(xy).
\end{equation}
And let us consider the nonlinear stochastic heat equation
%
\begin{equation}\label{eqmain}
\frac{\partial}{\partial t}u_t(x) = \frac\varkappa2 \,\frac{\partial^2}{
\partial x^2} u_t(x)+\sigma(u_t(x))\,\frac{\partial^2}{\partial t\,
\partial x}W(t,x),
\end{equation}
where $x\in\R$, $t>0$, $\sigma\dvtx\R\to\R$ is a nonrandom
and Lipschitz continuous function, $\varkappa> 0$
is a fixed viscosity parameter and the initial function $u_0\dvtx\R\to\R$
is bounded, nonrandom and measurable. The mixed partial derivative
$\partial^2 W(t,x)/(\partial t \,\partial x)$ is the so-called
``space--time white noise,'' and is defined
as a generalized Gaussian random field; see Chapter 2 of
Gelfand and Vilen\-kin~\cite{GV}, Section 2.7, for example.

It is well known that the stochastic heat equation
(\ref{eqmain}) has a (weak)
solution $\{u_t(x)\}_{t>0,x\in\R}$ that is jointly continuous;
it is also unique up to evanescence; see, for example, Chapter 3 of Walsh
\cite{Walsh}, (3.5), page 312.
And the solution can be written in mild form as
the (a.s.) solution to the following stochastic integral equation:
%
\begin{equation}\label{intromild}
u_t(x) = (p_t*u_0)(x) + \int_{(0,t)\times\R}
p_{t-s}(y-x)\sigma(u_s(y)) W(\d s \,\d y),
\end{equation}
where
%
\begin{equation}\label{p}
p_t(z) := \frac{1}{(2\pi\varkappa t)^{1/2}} \exp\biggl(-\frac
{z^2}{2\varkappa t}
\biggr)\qquad (t>0, z\in\R)
\end{equation}
denotes the free-space heat kernel, and
the final integral in (\ref{intromild})
is a stochastic integral in the
sense of Walsh~\cite{Walsh}, Chapter 2. Chapter 1 of
the minicourse by
Dalang et al.~\cite{Minicourse}
contains a quick introduction to the topic of stochastic PDEs
of the type considered here.

We are interested solely in the physically interesting case
that $u_0(x)\ge0$ for all $x\in\R$. In that case, a minor
variation of Mueller's comparison principle~\cite{Mueller}
implies that if in addition $\sigma(0)=0$,
then with probability one $u_t(x)\ge0$ for all $t>0$ and $x\in\R$;
see also Theorem 5.1 of Dalang et al.~\cite{Minicourse}, page 130,
as well as Theorem~\ref{thMueller} below.

We follow Foondun and Khoshnevisan~\cite{FK},
and say that the solution $u:=\{u_t(x)\}_{t>0,x\in\R}$ to (\ref{eqmain})
is (weakly) \textit{intermittent} if
%
\begin{equation}\label{interm}
0<\limsup_{t\to\infty} \frac1t\log\E(|u_t(x)|^\nu)
<\infty\qquad\mbox{for all $\nu\ge2$},
\end{equation}
where ``$\log$'' denotes the natural logarithm, to be concrete.
Here, we refer to property (\ref{interm}), if and when it holds,
as \textit{mathematical intermittency} [to be distinguished from
\textit{physical intermittency}, which is a phenomenological property
of an object that (\ref{eqmain}) is modeling].

If $\sigma(u)=\mathrm{const}\cdot u$ and $u_0$ is bounded
from above and below uniformly, then the work of
Bertini and Cancrini~\cite{BC} and
Mueller's comparison principle~\cite{Mueller} together
imply (\ref{interm}).
In the fully nonlinear case, Foondun and Khoshnevisan~\cite{FK}
discuss a connection to nonlinear renewal theory,
and use that connection to establish (\ref{interm})
under various conditions; for instance,
they have shown that (\ref{interm}) holds provided that
$\liminf_{|x|\to\infty}|\sigma(x)/x|>0$ and $\inf_{x \in\R}
u_0(x)$ is
sufficiently large.

If the $\limsup$ in (\ref{interm}) is a bona fide limit,
then we arrive at the usual description of intermittency
in the literature of mathematics and theoretical physics; see, for
instance,
Molchanov~\cite{Molch91} and Zeldovich et al.~\cite{ZMRS85,ZMRS88,ZRS}.

Mathematical intermittency is motivated strongly
by a vast physics literature on (physical) intermittency and
\textit{localization}, and many of\vadjust{\goodbreak} the references can be found in the combined
bibliographies of~\cite{FK,Molch91,ZMRS85,ZMRS88}.
Let us say a few more words about ``localization'' in the present
context.

It is generally accepted that if (\ref{interm}) holds, then
$u:=\{u_t(x)\}_{t>0,x\in\R}$ ought to undergo a separation of scales
(or ``pattern/period breaking'').
In fact, one can argue that property (\ref{interm})
implies that, as $t\to\infty$, the random function $x\mapsto u_t(x)$
starts to develop very tall peaks, distributed over small $x$-intervals
(see Section 2.4 of Bertini and Cancrini~\cite{BC} and the
Introduction of
the monograph by Carmona and Molchanov
\cite{CM}). This ``peaking property'' is called \textit{localization},
and is experienced with very high probability, provided that:
(i) the intermittency property (\ref{interm}) holds; and (ii) $t\gg1$.

Physical intermittency is expected to hold both in space and time, and
not only when $t\gg1$. And it is also expected that there are
physically-intermittent
processes, not unlike those studied in the present paper, which,
however, do not
satisfy the (mathematical) intermittency condition (\ref{interm}) on
Liapounov exponents; see, for example, the paper by
Cherktov et al.
\cite{CFKL95}.

Our wish is to better understand ``physical intermittency'' in the setting
of the stochastic heat equation (\ref{eqmain}). We are motivated strongly
by the literature on smooth finite-dimensional dynamical systems
(\cite{Ruelle}, Section 1.3), which ascribes
intermittency in part to ``chaos,'' or slightly more precisely,
sensitive dependence
on the initial state of the system.

In order to describe the contributions of this paper, we first recall
a consequence of a more general
theorem of Foondun and Khoshnevisan~\cite{FKAIHP}:
if $\sigma(0)=0$, and if $u_0$ is H\"older continuous of index
$>\frac12$ and has compact support, then for every $t>0$ fixed,
%
\begin{equation}\label{eqFK}
\limsup_{z\to\infty}\frac{1}{\log z}
\log\P\Bigl\{\sup_{x\in\R} u_t(x) > z\Bigr\}=-\infty.
\end{equation}
It follows in particular that the global maximum of the solution (at a
fixed time)
is a finite (nonnegative) random variable.

By contrast, one expects that if
%
\begin{equation}\label{u0bdd}
\inf_{x\in\R} u_0(x)>0,
\end{equation}
then the solution $u_t$ is unbounded for all $t>0$. Here we prove
that fact and a good deal more; namely, we demonstrate here that
there in fact exists a minimum rate of ``blowup'' that applies regardless
of the parameters of the problem.

A careful statement requires a
technical condition that turns out to be necessary as well
as sufficient.
In order to discover that condition, let
us consider the case that $u_0(x)\equiv\rho>0$
is a constant for all $x\in\R$.
Then, (\ref{u0bdd}) clearly holds;
but there can be no blowup if $\sigma(\rho)=0$.
Indeed, in that case the unique solution to the stochastic
heat equation is $u_t(x)\equiv\rho$, which is bounded.
Thus, in order to have an unbounded solution, we need at the
very least to consider the case that
$\sigma(x)\neq0$ for all $x > 0$.
[Note that $\sigma(0)=0$ is permitted.]
Instead, we will assume the following seemingly stronger, but
in fact more or less equivalent, condition
from now on:
%
\begin{equation}\label{sigmabdd}
\sigma(x)>0 \qquad\mbox{for all $x\in\R\setminus\{0\}$}.
\end{equation}

We are ready to present the first theorem of this paper.
Here and throughout we write
``$f(R)\succsim g(R)$ as $R\to\infty$'' in place of
the more cumbersome ``there exists a
nonrandom $C>0$ such that
$\liminf_{R\to\infty} f(R)/g(R)\ge C$.'' The largest such $C$
is called ``the constant in $\succsim$.''
We might sometimes also write ``$g(R)\precsim f(R)$''
in place of ``$f(R)\succsim g(R)$.'' And there is a corresponding
``constant in $\precsim$.''
%
\begin{theorem}\label{thmain}
Let $\{u_t(x)\}_{t>0,x\in\R}$ be a solution to
%
\begin{equation}\quad
\frac{\partial}{\partial t}u_t(x) = \frac\varkappa2 \,\frac{\partial^2}{
\partial x^2} u_t(x)+\sigma(u_t(x))\,\frac{\partial^2}{\partial t\,
\partial x}W(t,x)\qquad (t>0,x\in\R)
\end{equation}
written in mild form (\ref{intromild}), where the initial function
$u_0\dvtx\R\rightarrow\R$ is bounded, nonrandom and satisfies
%
\begin{equation}
\inf_{x\in\R} u_0(x)>0.
\end{equation}
Then, the following hold:
\begin{longlist}[(2)]
\item[(1)] If $\inf_{x\in\mathbf{R}}\sigma(x)\ge\varepsilon_0>0$
and $t>0$, then a.s.,
%
\begin{equation}
\sup_{x\in[-R,R]}u_t(x) \succsim
\frac{(\log R)^{1/6}}{\varkappa^{1/12}}\qquad
\mbox{as $R\to\infty$;}
\end{equation}
and the constant in $\succsim$ does not depend on $\varkappa$.
\item[(2)] If $\sigma(x)>0$ for all $x\in\R$ and there exists
$\gamma\in(0,1/6)$ such that
%
\begin{equation}\label{eqdecaysigma}
\lim_{|x|\to\infty} \sigma(x) \log(|x|)^{(1/6)-\gamma}
=\infty,
\end{equation}
then for all $t>0$
the following holds almost surely:
%
\begin{equation}
{\sup_{x\in[-R,R]}}|u_t(x)|
\succsim\frac{(\log R)^{\gamma}}{\varkappa^{1/12}}\qquad
\mbox{as $R\to\infty$;}
\end{equation}
and the constant in $\succsim$ does not depend on $\varkappa$.
\end{longlist}
\end{theorem}

Note in particular that if $\sigma$ is uniformly bounded below then a.s.,
%
\begin{equation}
\limsup_{|x|\to\infty}\frac{u_t(x)}{({\log}|x|)^{1/6}}\ge
\frac{\mathrm{const}}{\varkappa^{1/12}}.
\end{equation}
We believe that it is a somewhat significant fact that a rate $(\log
R)^{1/6}$ of blowup
exists that is valid for all $u_0$ and $\sigma$ in the first part of
Theorem~\ref{thmain}. However, the actual numerical estimate---that
is, the $(1/6)$th
power of the logarithm---appears to be
less significant, as the behavior in $\varkappa$
might suggest (see Remark~\ref{remkappa}). In fact,\vadjust{\goodbreak}
we believe that the actual blowup rate might depend critically on
the fine properties of the function $\sigma$. Next, we highlight this assertion
in one particularly interesting case. Here and throughout,
we write ``$f(R)\asymp g(R)$ as $R\to\infty$''
as shorthand for ``$f(R)\succsim g(R)$ and $g(R)\succsim
f(R)$ as $R\to\infty$.'' The two constants in the preceding
two $\succsim$'s are called the ``constants in $\asymp$.''
%
\begin{theorem}\label{thboundedsigma}
If $\sigma$ is uniformly bounded away from $0$ and $\infty$
and $t>0$, then
%
\begin{equation}
\sup_{x\in[-R,R]}u_t(x)
\asymp\frac{(\log R)^{1/2}}{\varkappa^{1/4}}\qquad
\mbox{a.s. as $R\to\infty$}.
\end{equation}
Moreover, for every fixed $\varkappa_0>0$,
the preceding constants in $\asymp$ do not depend on
$\varkappa\ge\varkappa_0$.
\end{theorem}

In particular, we find that if $\sigma$ is bounded uniformly away
from $0$ and $\infty$, then there exist constants
$c_*,c^*\in(0,\infty)$ such that
%
\begin{equation}
\frac{c_*}{\varkappa^{1/4}} \le
\limsup_{|x|\to\infty}\frac{u_t(x)}{({\log}|x|)^{1/2}}
\le\frac{c^*}{\varkappa^{1/4}}\qquad
\mbox{a.s.},
\end{equation}
uniformly for all $\varkappa\ge\varkappa_0$.

The preceding discusses the behavior in case $\sigma$ is
bounded uniformly away from $0$; that is, a uniformly-noisy stochastic heat
equation (\ref{eqmain}). In general, we can say little about the remaining
case that $\sigma(0)=0$. Nevertheless in the well-known
parabolic Anderson model, namely (\ref{eqmain}) with
$\sigma(x) = cx$ for some constant $c > 0$, we are able to
obtain some results (Theorem~\ref{thPAM})
that parallel Theorems~\ref{thmain} and~\ref{thboundedsigma}.
%
\begin{theorem} \label{thPAM}
If $\sigma(x) = cx$ for some $c>0$, then a.s.,
%
\begin{equation}
\log\sup_{x\in[-R,R]} u_t(x)
\asymp\frac{(\log R)^{2/3}}{\varkappa^{1/3}}\qquad
\mbox{as $R\to\infty$},
\end{equation}
and the constants in $\asymp$ do not depend on $\varkappa>0$.
\end{theorem}

Hence, when $\sigma(x)=cx$ we can find constants
$C_*,C^*\in(0,\infty)$ such that
\[
0<\limsup_{|x|\to\infty}\frac{u_t(x)}{\exp\{
C_*({\log}|x|)^{2/3}/\varkappa^{1/3}\}}
\le
\limsup_{|x|\to\infty}\frac{u_t(x)}{\exp\{
C^*({\log}|x|)^{2/3}/\varkappa^{1/3}\}}<\infty,
\]
almost surely.
%
\begin{remark}
Thanks to (\ref{intromild}), and since
Walsh stochastic integrals have zero mean,
it follows that $\E u_t(x)=(p_t*u_0)(x)$. In particular,
$\E u_t(x)\le\sup_{x \in\R} u_0(x)$ is
uniformly bounded. Since $u_t(x)$ is nonnegative,
it follows from Fatou's lemma that $\liminf_{|x|\to\infty}
u_t(x)<\infty$ a.s. Thus, the behavior described by
Theorem~\ref{thmain} is\vadjust{\goodbreak} one about the highly-oscillatory nature
of $x\mapsto u_t(x)$, valid for every fixed time $t>0$.
We will say a little more about this topic in
Appendix~\ref{secasinterm} below.
\end{remark}
%
\begin{remark} \label{remkappa}
We pay some attention to the powers of the viscosity parameter
$\varkappa$
in Theorems~\ref{thboundedsigma} and~\ref{thPAM}. Those
powers suggest that at least two distinct universality classes can be associated
to (\ref{eqmain}):
(i) when $\sigma$ is bounded uniformly away from zero and infinity,
the solution behaves as random walk in weakly-interacting
random environment; and
(ii) when $\sigma(x)=cx$ for some $c>0$, then the solution behaves
as objects that arise in some random matrix models.
\end{remark}
%
\begin{remark}
In~\cite{KPZ}, (2), Kardar, Parisi and
Zhang consider the solution $u$ to (\ref{eqmain}) and
apply formally the Hopf--Cole transformation
$u_t(x):=\exp(\lambda h_t(x))$ to deduce that
$\mathbf{h}:=\{h_t(x)\}_{t\ge0,x\in\R}$ satisfies the following\break
``SPDE'': for $t>0$ and $x\in\R$,
%
\begin{equation}\label{eqKPZ}
\frac{\partial}{\partial t} h_t(x) =
\frac{\varkappa}{2}\,\frac{\partial^2}{\partial x^2} h_t(x)
+\frac{\varkappa\lambda}{2}\biggl(
\frac{\partial}{\partial x}h_t(x)\biggr)^2 +
\frac{\partial^2}{\partial t \,\partial x}W(t,x).
\end{equation}
This is the celebrated ``KPZ equation,'' named after the authors of
\cite{KPZ}, and the random field $\mathbf{h}$
is believed to be a universal object (e.g., it is expected to arise
as a continuum limit of a large number of interacting particle systems).
Theorem~\ref{thPAM} implies that there exist positive and finite
constants $a_t$ and $A_t$---depending only on~$t$---such that
%
\begin{equation}
\frac{a_t}{\varkappa^{1/3}}
<\limsup_{|x|\to\infty} \frac{h_t(x)}{({\log}|x|)^{2/3}}<
\frac{A_t}{\varkappa^{1/3}}\qquad
\mbox{a.s. for all $t>0$}.
\end{equation}
This is purely formal, but only because the construction of $\mathbf{h}$
via $u$ is not rigorous. More significantly,
our proofs suggest strongly a kind of asymptotic
space--time scaling ``$|{\log x}|\approx
t^{\pm1/2}$.''
If so, then the preceding verifies that fluctuation exponent $1/z$
of $\mathbf{h}$ is $2/3$ under quite general conditions on the $h_0$. The latter
has been predicted by Kardar et al.~\cite{KPZ}, page 890, and proved
by Balazs, Quastel and Sepp\"al\"ainen~\cite{BQS} for a
special choice of $u_0$ (hence $h_0$) and $t \rightarrow0$.
\end{remark}

The proofs of our three theorems involve a fairly long series of technical
computations. Therefore, we conclude the \hyperref[intro]{Introduction} with a few
remarks on
the methods of proofs for the preceding three theorems in order to highlight
the ``pictures behind the proofs.''

Theorem~\ref{thmain} relies on two well-established techniques from
interacting particle systems~\cite{Durrett,Liggett}: namely, comparison
and coupling. Comparison reduces our problem to the case that $u_0$ is
a constant; at a technical level this uses Mueller's comparison
principle~\cite{Mueller}. And we use coupling on a few occasions:
first, we describe a two-step coupling of $\{u_t(x)\}_{t>0,x\in\R}$ to
the solution $\{v_t(x)\}_{t>0,x\in\R}$ of (\ref{eqmain})---using the\vadjust{\goodbreak}
same space--time white noise $\partial^2 W/(\partial t \,\partial
x)$---in the case that $\sigma$ is bounded below uniformly on $\R$. The
latter quantity [i.e., $\{v_t(x)\}_{t>0,x\in\R}$] turns out to be more
amenable to moment analysis than $\{u_t(x)\}_{t>0,x\in\R}$, and in this
way we obtain the following a priori estimate, valid for every $t>0$
fixed:
%
\begin{equation}\label{LBLB}
\log\inf_{x\in\R}
\P\{ u_t(x) \ge\lambda\} \succsim-\sqrt\varkappa
\lambda^6\qquad
\mbox{as $\lambda\to\infty$}.
\end{equation}
Theorem~\ref{thmain} follows immediately from this
and the Borel--Cantelli lemma,
provided that we prove that if $x$ and $x'$ are ``$O(1)$
distance apart,'' then $u_t(x)$ and $u_t(x')$ are ``approximately
independent.'' A quantitative version of this statement
follows from coupling $\{u_t(x)\}_{t>0,x\in\R}$ to
the solution $\{w_t(x)\}_{t>0,x\in\R}$
of a random evolution equation that can be
thought of as the ``localization'' of the original stochastic
heat equation (\ref{eqmain}). The localized approximation
$\{w_t(x)\}_{t>0,x\in\R}$ has the property that
$w_t(x)$ and $w_t(x')$ are (exactly) independent
for ``most'' values of $x$ and $x'$ that are $O(1)$
distance apart. And this turns out to be adequate for our needs.

Theorem~\ref{thboundedsigma} requires establishing
separately a lower and an upper bound on $\sup_{x\in[-R,R]}u_t(x)$.
Both bounds rely heavily on the following quantitative improvement of
(\ref{LBLB}): if $\sigma$ is bounded, then
%
\begin{equation}\label{LBUB}
\log\inf_{x\in\R}
\P\{ u_t(x) \ge\lambda\} \succsim-\sqrt\varkappa
\lambda^2\qquad
\mbox{as $\lambda\to\infty$}.
\end{equation}
And, as it turns out, the preceding lower bound will perforce imply
a corresponding upper estimate,
%
\begin{equation}\label{UBUB}
\log\inf_{x\in\R}
\P\{ u_t(x) \ge\lambda\} \precsim-\sqrt\varkappa
\lambda^2\qquad
\mbox{as $\lambda\to\infty$}.
\end{equation}
The derivation of the lower bound on $\sup_{x\in[-R,R]}u_t(x)$
follows closely the proof of Theorem~\ref{thmain},
after (\ref{LBUB}) and (\ref{UBUB}) are established.
Therefore, the remaining details will be omitted.

The upper bound on $\sup_{x\in[-R,R]}u_t(x)$ requires only
(\ref{UBUB}) and a well-known
quantitative version of the Kolmogorov continuity theorem.

Our proof of Theorem~\ref{thPAM} has a similar
flavor to that of Theorem~\ref{thmain}, for the lower bound, and
Theorem~\ref{thboundedsigma}, for the upper bound. We make
strong use of the moments formulas of Bertini and Cancrini
\cite{BC}, Theorem 2.6. [This is why we are only able to study the
linear equation in the case that $\sigma(0)=0$.]

Throughout this paper, we use the following abbreviation:
%
\begin{equation} \label{equstar}
u_t^*(R) := \sup_{x\in[-R,R]}u_t(x)\qquad
(R>0).
\end{equation}
We will also need the following elementary facts about the heat
kernel:
%
\begin{equation}\label{eqL2pt1}
\|p_s\|_{L^2(\R)}^2 =(4\pi\varkappa s)^{-1/2}\qquad
\mbox{for every $s>0$;}
\end{equation}
and therefore,
%
\begin{equation}\label{eqL2pt}
\int_0^t \|p_s\|_{L^2(\R)}^2 \,\d s =\sqrt{t/(\pi\varkappa)}\qquad
\mbox{for all $t\ge0$}.
\end{equation}
We will tacitly write $\mathrm{Lip}_\sigma$ for the optimal
Lipschitz constant of $\sigma$; that is,
%
\begin{equation}
\mathrm{Lip}_\sigma:= \sup_{-\infty<x\neq x'<\infty}
\biggl|\frac{\sigma(x)-\sigma(x')}{x-x'}\biggr|.
\end{equation}
Of course, $\mathrm{Lip}_\sigma$
is finite because $\sigma$ is Lipschitz continuous.
Finally, we use the following notation for the $L^\nu(\P)$
norm of a random variable $Z\in L^\nu(\P)$:
%
\begin{equation}
\|Z\|_\nu:= \{ \E( |Z|^\nu)\}^{1/\nu}.
\end{equation}

\section{Mueller's comparison principle and a reduction}\label{sec2}

Mueller's comparison principle~\cite{Mueller} is one of the
cornerstones of the theory of stochastic PDEs.
In its original form, Mueller's comparison principle is stated for
an equation that is similar to (\ref{eqmain}), but for two differences:
(i) $\sigma(z) :=\varkappa z$ for some $\varkappa>0$; and (ii) the variable
$x$ takes values in a compact interval such as $[0,1]$. In his Utah
Minicourse (\cite{Minicourse}, Theorem 5.1, page 130), Mueller
outlines how one can include also the more general functions
$\sigma$ of the type studied here. And in both cases, the proofs
assume that the initial
function $u_0$ has compact support. Below we state and
prove a small variation of the preceding
comparison principles that shows that Mueller's theory continues
to work when: (i) the variable $x$ takes values in $\R$; and (ii)
the initial function $u_0$ is not necessarily compactly supported.
%
\begin{theorem}[(Mueller's comparison principle)]\label{thMueller}
Let $u^{(1)}_0$ and $u^{(2)}_0$ denote two nonnegative bounded continuous
functions on $\R$ such that $u^{(1)}_0(x)\ge u^{(2)}_0(x)$ for all
$x\in\R$.
Let $u^{(1)}_t(x), u^{(2)}_t(x)$ be solutions
to (\ref{eqmain}) with respective initial functions $u^{(1)}_0$
and $u^{(2)}_0$. Then,
%
\begin{equation}
\P\bigl\{
u^{(1)}_t(x)\ge u^{(2)}_t(x)
\mbox{ for all $t>0$ and $x\in\R$} \bigr\}=1.
\end{equation}
\end{theorem}
\begin{pf}
Because the solution to (\ref{eqmain}) is continuous in $(t,x)$,
it suffices to prove that
%
\begin{equation}\label{eqcomp}
\P\bigl\{ u^{(1)}_t(x)\ge u^{(2)}_t(x)\bigr\}=1
\qquad\mbox{for all $t>0$ and $x\in\R$}.
\end{equation}
In the case that $u^{(1)}_0$ and $u^{(2)}_0$ both have bounded support,
the preceding is proved almost exactly as
in Theorem 3.1 of Mueller~\cite{Mueller}.
For general $u^{(1)}_0$ and $u^{(2)}_0$, we proceed as follows.

Let $v_0\dvtx\R\to\R_+$ be a bounded and measurable
initial function, and define a new initial function
$v_0^{[N]}\dvtx\R\to\R_+$ as
%
\begin{equation}
v_0^{[N]}(x) := \cases{
v_0(x), &\quad if $|x|\le N$,\cr
v_0(N)(-x+N+1), &\quad if $N<x<N+1$,\cr
v_0(-N)(x+N+1), &\quad if $-(N+1)<x<-N$,\cr
0, &\quad if $|x|\ge N+1$.}\vadjust{\goodbreak}
\end{equation}
Then, let $v_t^{[N]}(x)$ be the solution to (\ref{eqmain}) with
initial condition $v_0^{[N]}$. We claim that
%
\begin{equation}\label{eqtruncate}
\delta^{[N]}_t(x):= v_t(x) - v^{[N]}_t(x) \rightarrow0\qquad
\mbox{in probability as $N\rightarrow\infty$.}
\end{equation}
Let $u_t^{(1),[N]}$ and $u_t^{(2),[N]}$ denote the solutions to (\ref
{eqmain}) with initial conditions $u_0^{(1),[N]}$ and $u_0^{(2),[N]}$,
respectively, where the latter are defined similarly as $v_0^{[N]}$
above. Now, (\ref{eqtruncate}) has the desired result
because it shows that
$u^{(1),[N]}_t(x)\rightarrow u^{(1)}_t(x)$ and
$u^{(2),[N]}_t(x)\rightarrow u^{(2)}_t(x)$ in probability
as $N\rightarrow\infty$.
Since $u^{(1),[N]}_t(x)\ge u^{(2),[N]}_t(x)$ a.s. for all
$t>0$ and $x\in\R$, (\ref{eqcomp}) follows from taking limits.

In order to conclude we establish (\ref{eqtruncate});
in fact, we will prove that
%
\begin{equation}
\sup_{x\in\R}\sup_{t\in(0,T)}
\E\bigl(\bigl|\delta^{[N]}_t(x)\bigr|^2\bigr)=O(1/N)
\qquad\mbox{as $N\to\infty$}
\end{equation}
for all $T>0$ fixed. Recall that
%
\begin{eqnarray}\quad
\delta^{[N]}_t(x)&=&
\bigl( p_t*\delta^{[N]}_0\bigr)(x)\nonumber\\[-8pt]\\[-8pt]
&&{} +\int_{(0,t)\times\R}p_{t-s}(y-x)\bigl\{\sigma
( v_s(y) )-\sigma
\bigl( v_s^{[N]}(y) \bigr)\bigr\} W(\d s \,\d y).\nonumber
\end{eqnarray}
Because $( p_t*\delta^{[N]}_0 )(x)\le
2(\sup_{z\in\R} v_0(z)) \int_{|y|>N} p_t(y) \,\d y$,
a direct estimate of the latter stochastic integral yields
%
\begin{eqnarray}\label{ftm}\quad
\E\bigl( \bigl| \delta^{[N]}_t(x) \bigr|^2
\bigr) &\le& \mathrm{const}\cdot t^{-1/2}{\mathrm e}^{-N^2/(2t)}\nonumber\\
&&{} +\mathrm{const}\cdot\mathrm{Lip}_\sigma^2
\int_0^t\d s\int_{-\infty}^\infty\d y\,
p^2_{t-s}(y-x)\E\bigl(
\bigl| \delta^{[N]}_s(y)\bigr|^2
\bigr)\nonumber\\[-8pt]\\[-8pt]
&\le& \mathrm{const}\cdot t^{-1/2} {\mathrm e}^{-N^2/(2t)}\nonumber\\
&&{} +\mathrm{const}\cdot\mathrm{Lip}_\sigma^2
{\mathrm e}^{\beta t}M^{[N]}_t(\beta)
\int_0^\infty{\mathrm e}^{-\beta r}\| p_r\|_{L^2(\R)}^2 \,\d r,\nonumber
\end{eqnarray}
where $\beta>0$ is, for the moment, arbitrary and
%
\begin{equation}
M^{[N]}_t(\beta) := \sup_{s\in(0,t),y\in\R} \bigl[
{\mathrm e}^{-\beta s}\E\bigl(\bigl| v_s(y)-v_s^{[N]}(y)\bigr|^2\bigr)
\bigr].
\end{equation}
We multiply both sides of (\ref{ftm})
by $\exp(-\beta t)$ and take the supremum
over all $t\in(0,T)$ where $T>0$ is
fixed. An application of (\ref{eqL2pt1}) yields
%
\begin{equation}
M^{[N]}_T (\beta) \le\mathrm{const}\cdot\Bigl[
\sup_{t\in(0,T)}\bigl\{ t^{-1/2}{\mathrm e}^{-N^2/(2t)}\bigr\}
+ \beta^{-1/2} M^{[N]}_T (\beta) \Bigr].
\end{equation}
The quantity in $\sup_{t\in(0,T)}\{\cdots\}$ is proportional to
$1/N$ (with the constant of proportionality depending on $T$),
and the implied constant does not depend on $\beta$.\vspace*{1pt}
Therefore, it follows that if $\beta$ were selected sufficiently
large, then
$M^{[N]}_T (\beta) =O(1/N)$ as $N\to\infty$ for that choice
of $\beta$. This implies (\ref{eqtruncate}).
\end{pf}

Next we apply Mueller's comparison principle to make a
helpful simplification to our problem.

Because $\underline{B}:=\inf_{x \in\R} u_0(x)>0$ and
$\overline{B}:=\sup_{x \in\R} u_0(x)<\infty$, it follows from Theorem
\ref{thMueller} that almost surely,
%
\begin{equation}
\underline{u}_t(x)\le
u_t(x) \le\overline{u}_t(x)\qquad
\mbox{for all $t>0$ and $x\in\R$},
\end{equation}
where $\overline{u}$ solves the stochastic heat equation
(\ref{eqmain}) starting from initial function $\overline{u}_0(x):=
\overline{B}$, and $\underline{u}$ solves (\ref{eqmain})
starting from $\underline{u}_0(x):= \underline{B}$. This shows
that it suffices to prove Theorems~\ref{thmain},
\ref{thboundedsigma} and~\ref{thPAM} with $u_t(x)$
everywhere replaced by $\underline{u}_t(x)$ and
$\overline{u}_t(x)$. In other words, we can assume
without loss of generality that $u_0$ is identically a
constant.
In order to simplify the notation, we will
assume from now on that the mentioned constant is
one. A quick inspection of the ensuing proofs reveals
that this assumption is harmless.
Thus, from now on, we consider in place
of (\ref{eqmain}), the following parabolic stochastic PDE:
%
\begin{equation}\label{heat1}
\cases{\displaystyle
\frac{\partial}{\partial t} u_t(x)=\frac{\varkappa}{2}\,
\frac{\partial^2}{\partial x^2} u_t(x) + \sigma(u_t(x))
\,\frac{\partial^2}{\partial t\,\partial x}W(t,x), &\quad $(t > 0, x \in
\R)$,\vspace*{2pt}\cr
u_0(x)=1.}\hspace*{-32pt}
\end{equation}
We can write its solution in mild form as follows:
%
\begin{equation}\label{mild}
u_t(x) = 1 + \int_{[0,t]\times\R} p_{t-s}(y-x)
\sigma(u_s(y)) W(\d s\,\d y).
\end{equation}

\section{Tail probability estimates} \label{sectailestimates}

In this section we derive the following corollary
which estimates the tails of the distribution
of $u_t(x)$, where $u_t(x)$ solves (\ref{heat1}) and (\ref{mild}). In fact,
Corollary~\ref{comomentsU}, Propositions~\ref{prlowertail} and
\ref{prsigmadecay} below readily imply the following:
%
\begin{corollary}\label{cobds}
If $\inf_{z\in\R}[ \sigma(z)]>0$, then for all $t>0$,
%
\begin{equation}\label{eqtails1}
-\sqrt\varkappa\lambda^{6} \precsim
\log\P\{ |u_t(x)|\ge\lambda\}\precsim
-\sqrt\varkappa(\log\lambda)^{3/2},
\end{equation}
uniformly for $x\in\R$ and $\lambda>{\mathrm e}$. And the constants
in $\precsim$ do not depend on $\varkappa$.

If (\ref{eqdecaysigma})
holds for some $\gamma\in(0,1/6)$, then for all $t>0$,
%
\begin{equation}\label{eqtails}
-\varkappa^{1/(12\gamma)} \lambda^{1/\gamma}\precsim
\log\P\{ |u_t(x)|\ge\lambda\} \precsim
-\sqrt\varkappa(\log\lambda)^{3/2},
\end{equation}
uniformly for all $x\in\R$ and $\lambda>{\mathrm e}$. And the constants
in $\precsim$ do not depend on $\varkappa$.
\end{corollary}

\subsection{An upper-tail estimate} \label{secuppertail}
We begin by working toward the upper
bound in Corollary~\ref{cobds}.
%
\begin{lemma}\label{lemmomentsU}
Fix $T>0$, and define
$a:=T(\mathrm{Lip}_\sigma\vee1)^4 / (2\varkappa)$.
Then, for all real numbers $k\ge1$,
\[
\sup_{x\in\R}\sup_{t\in[0,T]}
\E( |u_t(x)|^k) \le C^k {\mathrm e}^{ak^3}\qquad
\mbox{where }
C:=8\biggl( 1+\frac{|\sigma(0)|}{
2^{1/4}(\mathrm{Lip}_\sigma\vee1)}\biggr).
\]
\end{lemma}
\begin{pf}
We follow closely the proof of Theorem 2.1 of
\cite{FK}, but matters simplify considerably in the present, more
specialized, setting.

First of all, we note that because $u_0\equiv1$ is spatially homogeneous,
the distribution of $u_t(x)$
does not depend on $x$; this property was observed earlier
by Dalang~\cite{Dalang99}, for example.

Therefore, an application of
Burkholder's inequality, using the Carlen--Kree bound~\cite{CK}
on Davis's optimal constant~\cite{Davis} in the Burkholder--Davis--Gundy
inequality~\cite{Burkholder,BDG,BG} and Minkowski's inequality imply
the following: for all $t\ge0$,
$\beta>0$ and $x\in\R$,
%
\begin{eqnarray}
\| u_t(x)\|_k
&\le&1 +\biggl\| \int_{[0,t]\times\R} p_{t-s}(y-x)
\sigma(u_s(y)) W(\d s\,\d y)\biggr\|_k\nonumber\\
&\le&1+ {\mathrm e}^{\beta t}2\sqrt{k}\Bigl( |\sigma(0)|+\mathrm{Lip}_\sigma
\sup_{r\ge0}[ {\mathrm e}^{-2\beta r}\|u_r(x)\|_k ]
\Bigr)\nonumber\\[-8pt]\\[-8pt]
&&\hspace*{15.6pt}{}\times\biggl(\int_{0}^{\infty} {\mathrm e}^{-2\beta s} \|p_s\|^2_2 \,\d s
\biggr)^{1/2} \nonumber\\
&=&1+\frac{\sqrt{k}{\mathrm e}^{\beta t}}{(8\varkappa\beta)^{1/4}}
\Bigl( |\sigma(0)|+\mathrm{Lip}_\sigma
\sup_{r\ge0}[ {\mathrm e}^{-2\beta r}\|u_r(x)\|_k ] \Bigr).\nonumber
\end{eqnarray}
%
See Foondun and Khoshnevisan~\cite{FK}, Lemma 3.3, for
the details of the derivation of such an estimate.
(Although Lemma 3.3 of~\cite{FK} is stated for even integers $k\ge2$,
a simple variation on the proof of that lemma implies the result
for general $k\ge1$; see Conus and Khoshnevisan~\cite{CKh}.)
It follows that
%
\begin{equation}
\psi(\beta,k):=\sup_{t\ge0}[ {\mathrm e}^{-\beta t}\|u_t(x)\|
_k]
\end{equation}
satisfies
%
\begin{equation}
\psi(\beta,k) \le1 + \frac{\sqrt{k}}{(4\varkappa\beta)^{1/4}}
\bigl( |\sigma(0)|+\mathrm{Lip}_\sigma\psi(\beta,k) \bigr).
\end{equation}
If\vspace*{1pt} $\mathrm{Lip}_\sigma=0$, then clearly $\psi(\beta,k)<\infty
$. If
$\mathrm{Lip}_\sigma>0$, then $\psi(\beta,k)<\infty$ for all
$\beta> k^2\operatorname{Lip}_\sigma^4/(4\varkappa)$;
therefore, the preceding proves that if
$\beta> k^2\operatorname{Lip}_\sigma^4/(4\varkappa)$,
then
%
\begin{equation}
\psi(\beta,k) \le\frac{1}{1-(\sqrt{k}
\operatorname{Lip}_\sigma/(4\varkappa\beta)^{1/4})}
\cdot\biggl( 1+ \frac{\sqrt{k}|\sigma(0)|}{(4\varkappa\beta
)^{1/4}}\biggr).
\end{equation}
We apply this with $\beta:=k^2(\mathrm{Lip}_\sigma\vee1)^4 /
(2\varkappa)$ to obtain the lemma.
\end{pf}
%
\begin{remark}
In the preceding results, the term $\mathrm{Lip}_\sigma\vee1$
appears in place of the more natural quantity $\mathrm{Lip}_\sigma$
only because it can happen that $\mathrm{Lip}_\sigma= 0$.
In the latter case, $\sigma$ is a constant function,
and the machinery of Lemma~\ref{lemmomentsU} is not needed
since $u_t(x)$ is a centered Gaussian process with a variance that can be
estimated readily. (We remind the reader that the case where $\sigma$ is
a constant is covered
by Theorem~\ref{thboundedsigma}; see Section~\ref{secboundedsigma}.)
\end{remark}

Next we describe a real-variable lemma that shows how to transform the
moment estimate of Lemma~\ref{lemmomentsU} into subexponential
moment estimates.
%
\begin{lemma}\label{lemexplog}
Suppose $X$ is a nonnegative random variable that satisfies
the following: there exist finite numbers
$a,C>0$ and $b>1$ such that
%
\begin{equation}
\E(X^k) \le C^k {\mathrm e}^{a k^b} \qquad\mbox{for all real numbers
$k\ge1$}.
\end{equation}
Then, $\E\exp\{\alpha(\log_+ X)^{b/(b-1)} \}<\infty$---for
$\log_+u := \log(u\vee{\mathrm e})$---provided that
%
\begin{equation}
0<\alpha< \frac{1- b^{-1}}{(ab)^{1/(b-1)}}.
\end{equation}
\end{lemma}

Lemmas~\ref{lemmomentsU},~\ref{lemexplog} and Chebyshev's
inequality together
imply the following result.
%
\begin{corollary}\label{comomentsU}
Choose and fix $T>0$, and define
$c_0 :=\sqrt{2/3}\approx0.8165$.
Then for all $\alpha<c_0\sqrt{\varkappa}/(\sqrt{T}(\mathrm{Lip}_\sigma\vee1))$,
%
\begin{equation}
\sup_{x\in\R}\sup_{t\in[0,T]}
\E\bigl({\mathrm e}^{\alpha(\log_+ u_t(x))^{3/2}}\bigr)<\infty.
\end{equation}
Consequently,
%
\begin{equation}\quad
\limsup_{\lambda\uparrow\infty}
\frac{1}{(\log\lambda)^{3/2}}
\sup_{x\in\R}\sup_{t\in[0,T]}\log
\P\{ u_t(x)>\lambda\} \le
\frac{-c_0 \sqrt{\varkappa}}{\sqrt{T}(\mathrm{Lip}_\sigma\vee1)}.
\end{equation}
\end{corollary}

We skip the derivation of Corollary~\ref{comomentsU} from
Lemma~\ref{lemexplog}, as it is immediate. The result holds uniformly
in $t \in[0,T]$
and $x \in\R$ as the constants $a$ and $C$ in Lem\-ma~\ref
{lemmomentsU} are independent of
$t$ and $x$. Instead we verify Lemma~\ref{lemexplog}.
\begin{pf*}{Proof of Lemma~\ref{lemexplog}}
Because
%
\begin{equation}
\biggl[ \log_+\biggl(\frac{X}{C}\biggr) \biggr]^{b/(b-1)} \le
2^{b/(b-1)}\cdot\bigl\{
(\log_+ X)^{b/(b-1)} + (\log_+ C)^{b/(b-1)}\bigr\},\hspace*{-35pt}\vadjust{\goodbreak}
\end{equation}
we can assume without loss of generality that $C=1$; for otherwise
we may consider $X/C$ in place of $X$ from here on.

For all $z>{\mathrm e}$, Chebyshev's inequality implies that
%
\begin{equation}
\P\bigl\{ {\mathrm e}^{\alpha(\log_+ X)^{b/(b-1)}} > z\bigr\}\le
{\mathrm e}^{-\max_k g(k)},
\end{equation}
where
%
\begin{equation}
g(k) := k\biggl(\frac{\log z}{\alpha}\biggr)^{(b-1)/b} - ak^b.
\end{equation}
One can check directly that $\max_k g(k) =c\log z$,
where
%
\begin{equation}
c:= \frac{1- b^{-1}}{\alpha\cdot(ab)^{1/(b-1)}}.
\end{equation}
Thus, it follows that $\P\{ \exp[\alpha(\log_+ X)^{b/(b-1)}]>z\}=O(z^{-c})$
for $z\to\infty$. Consequently,
$\E\exp\{\alpha(\log_+ X)^{b/(b-1)}\}<\infty$ as long as $c>1$;
this is equivalent to the statement of the lemma.
\end{pf*}

\subsection{Lower-tail estimates}
In this section we proceed to estimate the tail of
the distribution of $u_t(x)$ from below. 
We first consider the simplest case in which $\sigma$
is bounded uniformly from below, away from zero.
%
\begin{proposition}\label{prlowermoments}
If $\varepsilon_0:=\inf_{z\in\mathbf{R}}\sigma(z)>0$, then for all $t>0$,
%
\begin{equation}
\inf_{x\in\R}\E( | u_t(x)|^{2k}) \ge
\bigl( \sqrt{2}+o(1) \bigr) (\mu_tk)^k\qquad
(\mbox{as $k\to\infty$})
\end{equation}
where the ``$o(1)$'' term depends only on $k$, and
%
\begin{equation} \label{eqndefmu}
\mu_t := \frac{2}{{\mathrm e}}\cdot\varepsilon_0^2\sqrt{\frac{t}{\pi
\varkappa}}.
\end{equation}
\end{proposition}
\begin{pf}
Because the initial function in (\ref{heat1}) is $u_0(x)\equiv1$,
it follows that the distribution of $u_t(x)$ does not depend on
$x$; see Dalang~\cite{Dalang99}. Therefore, the ``$\inf$'' in the
statement of
the proposition is superfluous.

Throughout, let us fix $x\in\R$ and $t>0$. Now we may consider
a mean-one martingale $\{M_\tau\}_{0\le\tau\le t}$ defined as follows:
%
\begin{equation}\label{eqMtau}\qquad
M_\tau:= 1 + \int_{(0,\tau)\times\R} p_{t-s}(y-x)\sigma(u_s(y))
W(\d s \,\d y)\qquad (0\le\tau\le t).
\end{equation}
The quadratic variation of this martingale is
%
\begin{equation}
\langle M\rangle_\tau=\int_0^\tau\d s\int_{-\infty}^\infty\d y\,
p_{t-s}^2(y-x)\sigma^2(u_s(y))\qquad
(0\le\tau\le t).
\end{equation}
Therefore, by It\^o's formula, for all positive integers $k$,
and for every $\tau\in[0,t]$,
%
\begin{eqnarray}\label{eqIto}
M^{2k}_\tau
&=& 1+2k\int_0^\tau M_s^{2k-1} \,\d M_s
+\pmatrix{2k\cr2}\int_0^\tau M_s^{2(k-1)} \,\d\langle M\rangle
_s\nonumber\\
&=&1+2k\int_0^\tau M_s^{2k-1} \,\d M_s\\
&&{} +\pmatrix{2k\cr2}
\int_0^\tau M_s^{2(k-1)}\,\d s\int_{-\infty}^\infty\d y\,
p_{t-s}^2(y-x)\sigma^2(u_s(y)).\nonumber
\end{eqnarray}
%
By the assumption of the lemma, $\sigma(u_s(y)) \ge\varepsilon_0$ a.s.
Therefore,
%
\begin{eqnarray}\label{eqItoLB}\qquad
M^{2k}_\tau&\ge&1+2k\int_0^\tau M_s^{2k-1} \,\d M_s
+\pmatrix{2k\cr2}\varepsilon_0^2\int_0^\tau M_s^{2(k-1)}
\cdot\|p_{t-s}\|_{L^2(\R)}^2 \,\d s\nonumber\\[-8pt]\\[-8pt]
&=& 1+2k\int_0^\tau M_s^{2k-1} \,\d M_s +
\pmatrix{2k\cr2} \varepsilon_0^2 \int_0^\tau
\frac{M_s^{2(k-1)}}{(4\pi\varkappa(t-s))^{1/2}} \,\d s.\nonumber
\end{eqnarray}
We set $\tau:=t$ and then
take expectations to find that
%
\begin{eqnarray}\label{eqrecineq}
\E( M^{2k}_t) &\ge& 1+\pmatrix{2k\cr2}\varepsilon_0^2
\int_0^t \E\bigl(M^{2(k-1)}_s\bigr)\,\frac{\d s}{(4\pi\varkappa
(t-s))^{1/2}}\nonumber\\[-8pt]\\[-8pt]
&=& 1+\pmatrix{2k\cr2}\varepsilon^2_0\int_0^t
\E\bigl(M^{2(k-1)}_s\bigr) \nu(t,\d s),\nonumber
\end{eqnarray}
where the measures $\{\nu(t,\cdot)\}_{t>0}$ are defined as
%
\begin{equation}
\nu(t,\d s):= \frac{\one_{(0,t)}(s)}{(4\pi\varkappa
(t-s))^{1/2}}
\,\d s.
\end{equation}
We may iterate the preceding in order to obtain
%
\begin{eqnarray}\label{eqMLB}
\E( M^{2k}_t )
&\ge& 1 + \sum_{l=0}^{k-1} a_{l,k}
\varepsilon_0^{2(l+1)}\cdot
\int_0^t \nu(t,\d s_1) \int_0^{s_1}
\nu(s_1,\d s_2)
\cdots\nonumber\\[-8pt]\\[-8pt]
&&\hspace*{33pt}{}\times\int_0^{s_l}
\nu(s_l,\d s_{l+1}),\nonumber
\end{eqnarray}
where
%
\begin{equation}
a_{l,k} := \prod_{j=0}^l \pmatrix{2k-2j\cr2}\qquad
\mbox{for $0\le l<k$}
\end{equation}
and $s_0 :=t$. The right-hand side of (\ref{eqMLB}) is exactly equal
to $\E(M_t^{2k})$
in the case where $\sigma(z)\equiv\varepsilon_0$ for all $z\in\R$. Indeed,
the same computation as above works with identities all the way through.
In other words,
%
\begin{equation}
\E(|u_t(x)|^{2k})=
\E( M_t^{2k}) \ge\E( \eta_t(x)^{2k}),
\end{equation}
where
%
\begin{equation}
\eta_t(x) := 1+\varepsilon_0\cdot\int_{(0,t)\times\R} p_{t-s}(y-x)
W(\d s \,\d y).
\end{equation}
We define
%
\begin{equation}
\zeta_t(x):=\varepsilon_0\cdot\int_{(0,t)\times\R} p_{t-s}(y-x)
W(\d s\,\d y),
\end{equation}
so that $\eta_{t}(x)=1+\zeta_t(x)$. Clearly,
%
\begin{equation}
\E(\eta_t(x)^{2k})\ge\E(\zeta_t(x)^{2k}).
\end{equation}
Since $\zeta$ is a centered Gaussian process,
%
\begin{equation}
\E(\zeta_t(x)^{2k}) = [\E(\zeta
_t(x)^2
)]^{k}
\cdot\frac{(2k)!}{k!\cdot2^{k}}
\end{equation}
and
%
\begin{equation}
\E(\zeta_t(x)^2) = \varepsilon_0^2\cdot
\int_0^t\d s\int_{-\infty}^\infty\d y\,
p_{t-s}^2(y-x) = \varepsilon_0^2\cdot\sqrt{\frac{t}{\pi\varkappa}};
\end{equation}
see (\ref{eqL2pt}). The proposition follows from these
observations and Stirling's formula.
\end{pf}

We can now use Proposition~\ref{prlowermoments} to obtain a lower
estimate on the tail of the distribution of $u_t(x)$.
%
\begin{proposition}\label{prlowertail}
If there exists $\varepsilon_0>0$
such that $\sigma(x)\ge\varepsilon_0$ for all $x\in\mathbf{R}$, then
there exists a universal constant $C\in(0,\infty)$ such that for
all $t>0$,
%
\begin{equation}
\liminf_{\lambda\to\infty}\frac{1}{\lambda^{6}}
\inf_{x\in\R}
\log\P\{|u_t(x)|\ge\lambda\}
\ge- C \frac{(\mathrm{Lip}_\sigma\vee1)^4 \sqrt{\varkappa
}}{\varepsilon_0^6 \sqrt{t}}.
\end{equation}
\end{proposition}
\begin{pf}
Choose and fix $t>0$ and $x\in\R$.
We apply the celebrated
Paley--Zygmund inequality in the following form:
for every integer $k\ge1$,
%
\begin{eqnarray}
\E(|u_t(x)|^{2k})
&\le&\E\bigl(|u_t(x)|^{2k} ;
|u_t(x)|\ge\tfrac12 \|u_t(x)\|_{2k}\bigr)+
\tfrac12 \E(|u_t(x)|^{2k})\nonumber\hspace*{-35pt}\\[-8pt]\\[-8pt]
&\le&\sqrt{\E(|u_t(x)|^{4k})\P\bigl\{
|u_t(x)|\ge\tfrac12 \|u_t(x)\|_{2k}\bigr\}}+
\tfrac12 \E(|u_t(x)|^{2k}).\nonumber\hspace*{-35pt}
\end{eqnarray}
This yields the following bound:
%
\begin{eqnarray}\label{eqpzbound}\qquad
\P\biggl\{
|u_t(x)|\ge\frac12\|u_t(x)\|_{2k}\biggr\}
&\ge&\frac{[
\E(|u_t(x)|^{2k})]^2}{%
4\E(|u_t(x)|^{4k})}\nonumber\\[-8pt]\\[-8pt]
&\ge&\exp\biggl(- \frac{64 t(\mathrm{Lip}_\sigma\vee
1)^4}{\varkappa
} k^3\bigl(1+o(1)\bigr)
\biggr)\nonumber
\end{eqnarray}
as $k\to\infty$; see Lemma~\ref{lemmomentsU} and
Proposition~\ref{prlowermoments}. 
Another application of Proposition~\ref{prlowermoments} 
shows that $\|u_t(x)\|_{2k}
\ge(1+o(1))(\mu_t k)^{1/2}$ as $k\to\infty$ where $\mu
_t$ is
defined in (\ref{eqndefmu}). 
This implies as $k\to\infty$,
%
\begin{equation} \label{eq-128}\quad
\P\biggl\{ |u_t(x)| \ge\frac{1}{2} (\mu_t k)^{1/2} \biggr\}
\ge\exp\biggl[-\frac{64 t(\mathrm{Lip}_\sigma\vee1)^4}{
\varkappa} k^3\bigl(1+o(1)\bigr)
\biggr].
\end{equation}
The proposition follows from this by setting $k$
to be the smallest possible integer that
satisfies $(\mu_t k)^{1/2} \ge\lambda$.
\end{pf}

Now, we study the tails of the distribution of $u_t(x)$
under the conditions of part~(2) of Theorem~\ref{thmain}.
%
\begin{proposition}\label{prsigmadecay}
Suppose $\sigma(x)>0$ for all $x\in\mathbf{R}$ and
(\ref{eqdecaysigma}) holds
for some $\gamma\in(0,1/6)$. Then
%
\begin{equation}\quad
\liminf_{\lambda\to\infty} 
\inf_{x\in\mathbf{R}}\frac{\log\P\{
|u_t(x)| >\lambda\}}{\lambda^{1/\gamma}}
\ge- C \biggl( \frac{(\mathrm{Lip}_\sigma\vee1)^{2/3} \varkappa
^{1/12}} {t^{1/12}}
\biggr)^{1/\gamma},
\end{equation}
where $C\in(0,\infty)$
is a constant that depends only on $\gamma$.
\end{proposition}
\begin{pf}
For every integer $N\ge1$,
define
%
\begin{equation}
\sigma^{(N)}(x) := \cases{
\sigma(x), &\quad if $|x|\le N$,\cr
\sigma(-N), &\quad if $x<-N$,\cr
\sigma(N), &\quad if $x>N$.}
\end{equation}
It can be checked directly that $\sigma^{(N)} $ is a Lipschitz
function, and that in fact: (i)~$\mathrm{Lip}_{\sigma^{(N)}}
\le\mathrm{Lip}_{\sigma}$; and (ii) and
$\inf_{z\in\mathbf{R}}\sigma^{(N)}(z)>0$.

Let $u^{(N)}_t(x)$ denote the solution to (\ref{heat1}), when $\sigma$
is replaced by $\sigma^{(N)}$. We first establish
the bound
%
\begin{equation}
\E\bigl(\bigl| u^{(N)}_t(x)-u_t(x)\bigr|^2\bigr)=
O(N^{-2})
\qquad\mbox{as $N\to\infty$}.
\end{equation}

Let us observe, using the mild representation of the solution
to (\ref{heat1}), that
%
\begin{equation}\label{equ-uN}
\E\bigl(\bigl| u^{(N)}_t(x)-u_t(x)\bigr|^2\bigr)
\le2(\mathcal{T}_1+\mathcal{T}_2),
\end{equation}
where
%
\begin{eqnarray}
\mathcal{T}_1 :\!&=& \E\biggl( \biggl|
\int_{(0,t)\times\R}p_{t-s}(y-x)
\bigl[ \sigma^{(N)}\bigl(u^{(N)}_s(y)\bigr)-\sigma
\bigl( u^{(N)}_s(y)\bigr) \bigr] W(\d s\,\d y)
\biggr|^2\biggr)\hspace*{-32pt}\nonumber\\
&=& \int_0^t\d s\int_{-\infty}^\infty\d y\, p_{t-s}^2(y-x)
\E\bigl(\bigl| \sigma^{(N)}\bigl(u^{(N)}_s(y)\bigr)-\sigma
\bigl( u^{(N)}_s(y)\bigr) \bigr|^2\bigr)\quad\mbox{and}
\hspace*{-32pt}\nonumber\\[-8pt]\\[-8pt]
\mathcal{T}_2 :\!&=& \E\biggl( \biggl|
\int_{(0,t)\times\R}p_{t-s}(y-x) \bigl[
\sigma\bigl( u^{(N)}_s(y) \bigr)
-\sigma( u_s(y) ) \bigr] W(\d s \,\d y)\biggr|^2
\biggr)\hspace*{-32pt}\nonumber\\
&=& \int_0^t\d s\int_{-\infty}^\infty\d y\, p_{t-s}^2(y-x) \E
\bigl(\bigl|
\sigma\bigl( u^{(N)}_s(y) \bigr)
-\sigma( u_s(y) ) \bigr|^2\bigr).\hspace*{-32pt}\nonumber
\end{eqnarray}

We can estimate the integrand of $\mathcal{T}_1$ by the following:
%
\begin{eqnarray}\label{equN-uN}
&&\E\bigl(\bigl| \sigma^{(N)}\bigl(u^{(N)}_s(y)\bigr)-\sigma
\bigl( u^{(N)}_s(y)\bigr) \bigr|^2\bigr)\nonumber\\
&&\qquad\le\mathrm{Lip}_{\sigma}^2\cdot\E\bigl(
\bigl| N-u^{(N)}_s(y)\bigr|^2;
u^{(N)}_s(y)>N\bigr)\nonumber\\[-8pt]\\[-8pt]
&&\qquad\quad{} + \mathrm{Lip}_{\sigma}^2\cdot\E\bigl(
\bigl| -N-u^{(N)}_s(y)\bigr|^2; u^{(N)}_s(y)<-N\bigr)\nonumber\\
&&\qquad\le4\operatorname{Lip}_\sigma^2\cdot\E\bigl( \bigl|
u^{(N)}_s(y)\bigr|^2
; \bigl| u^{(N)}_s(y) \bigr| >N\bigr).\nonumber
\end{eqnarray}
We first apply the Cauchy--Schwarz inequality and then Chebyshev's
inequality (in this order) to conclude that
%
\begin{eqnarray}\label{equN-uN2}
\E\bigl(\bigl| \sigma^{(N)}\bigl(u^{(N)}_s(y)\bigr)-\sigma
\bigl( u^{(N)}_s(y)\bigr) \bigr|^2\bigr)
&\le&\frac{4\operatorname{Lip}_\sigma^2}{N^2}
\cdot\E\bigl( \bigl| u^{(N)}_s(y)\bigr|^4\bigr)\nonumber\\[-8pt]\\[-8pt]
&=& O(N^{-2}) \qquad\mbox{as $N\to\infty$},\nonumber
\end{eqnarray}
uniformly for all $y\in\R$ and $s\in(0,t)$. Indeed, Lemma \ref
{lemmomentsU}
ensures that $\E[|u_s^{(N)}(y)|^4]$ is bounded in $N$, because
$\lim_{N\to\infty}
\mathrm{Lip}_{\sigma^{(N)}}= \mathrm{Lip}_{\sigma}$.
This implies readily that
$\mathcal{T}_1=O(N^{-2})$ as $N\to\infty$.

Next we turn to $\mathcal{T}_2$; thanks to
(\ref{eqL2pt}), the quantity $\mathcal{T}_2$
can be estimated as follows:
%
\begin{eqnarray}\label{eqT2est}
\mathcal{T}_2 &\le&
\mathrm{Lip}_\sigma^2\cdot\int_0^t\d s\int_{-\infty}^\infty\d y\,
p_{t-s}^2(y-x) \E\bigl(\bigl|
u^{(N)}_s(y)-u_s(y)\bigr|^2\bigr)\nonumber\\[-8pt]\\[-8pt]
&\le&\mathrm{const}\cdot\int_0^t
\frac{\mathcal{M}(s)}{\sqrt{t-s}},\nonumber
\end{eqnarray}
where
%
\begin{equation}
\mathcal{M}(s) := \sup_{y\in\R}\E\bigl(
\bigl| u^{(N)}_s(y)-u_s(y)\bigr|^2\bigr)\qquad
(0\le s\le t).
\end{equation}
Notice that the implied constant in (\ref{eqT2est}) does not depend
on $t$.

We now combine our estimates for $\mathcal{T}_1$
and $\mathcal{T}_2$ to conclude that
%
\begin{eqnarray}
\mathcal{M}(s) &\le& \frac{\mathrm{const}}{N^2} + \mathrm{const} \cdot
\int_0^s\frac{\mathcal{M}(r)}{
\sqrt{s-r}} \,\d r \qquad(0\le s\le t) \nonumber\\[-8pt]\\[-8pt]
&\le& \mathrm{const}\cdot\biggl\{
\biggl( \int_0^s [\mathcal{M}(r)]^{3/2} \,\d r\biggr)^{2/3} + \frac{1}{N^2}
\biggr\},\nonumber
\end{eqnarray}
thanks to H\"older's inequality. We emphasize that the implied constant
depends only on the Lipschitz constant of $\sigma$, the variable $t$
and the diffusion constant $\varkappa$.
Therefore,
%
\begin{equation}
[\mathcal{M}(s)]^{3/2} \le\mathrm{const} \cdot\biggl\{
\int_0^s[\mathcal{M}(r)]^{3/2} \,\d r + \frac{1}{N^3} \biggr\},
\end{equation}
uniformly for $s\in(0,t)$. Gronwall's inequality
then implies the bound $\mathcal{M}(t) =O(N^{-2})$, valid
as $N\to\infty$.

Now we proceed with the proof of Proposition~\ref{prsigmadecay}.
For all $N\ge1$, the function $\sigma^{(N)}$ is bounded below. Let
$\varepsilon(N)$ be such that $\sigma^{(N)}(x)\ge\varepsilon(N)$ for all
$x \in\R$. Let $D:=D_t:=(4t/({\mathrm e}^2\pi\varkappa))^{1/4}$.
According to
the proof of Proposition~\ref{prlowertail}, specifically (\ref{eq-128})
applied to $u^{(N)}$, we have
%
\begin{eqnarray}\label{equngamma}
\P\biggl\{ |u_t(x)| \ge\frac{D}{4}\varepsilon(N)k^{1/2} \biggr\}
&\ge&\exp\biggl[- \frac{64 t(\mathrm{Lip}_\sigma\vee1)^4}{
\varkappa} k^3\bigl(1+o(1)\bigr)\biggr]\nonumber\hspace*{-35pt}\\[-8pt]\\[-8pt]
&&{} -\P\biggl\{\bigl\vert u_t(x)-u^{(N)}_t(x)
\bigr\vert\ge\frac{D}{4}\varepsilon(N)k^{1/2}\biggr\}.
\nonumber\hspace*{-35pt}
\end{eqnarray}
Thanks to (\ref{eqdecaysigma}), we can write
%
\begin{equation}
\varepsilon(N) \gg(\log N)^{-(1/6-\gamma)}\qquad
\mbox{as $N\to\infty$},
\end{equation}
using standard notation. Therefore, if
we choose
%
\begin{equation}
N :=\biggl\lfloor
\exp\biggl\{\frac{64 t(\mathrm{Lip}_\sigma\vee1)^4 k^3
}{\varkappa
}\biggr\}
\biggr\rfloor,
\end{equation}
then we are led to the bound
%
\begin{equation}
\varepsilon(N) \gg\biggl(\frac{64 t (\mathrm{Lip}_\sigma\vee1)^4}{
\varkappa}\biggr)^{-(1/6-\gamma)} k^{3\gamma-(1/2)}.
\end{equation}
We can use Chebyshev's inequality in order to estimate
the second term on the right-hand side of (\ref{equngamma}).
In this way we obtain the following:
%
\begin{eqnarray}
\P\biggl\{ |u_t(x)| \ge\frac{\tilde{D}}{4}k^{3\gamma} \biggr\}
&\ge&\exp\biggl\{-\frac{64 t(\mathrm{Lip}_\sigma\vee1)^4}{
\varkappa} k^3\bigl(1+o(1)\bigr) \biggr\}\nonumber\\[-8pt]\\[-8pt]
&&{}-\frac{1}{C_1N^2k^{6\gamma}},\nonumber
\end{eqnarray}
where
%
\begin{equation}
\tilde{D} := D \biggl\{\frac{\varkappa}{
64 t (\mathrm{Lip}_\sigma\vee1)^4}\biggr\}^{(1/6)-\gamma},
\end{equation}
and $C_1$ is a constant that depends only on $t$,
$\mathrm{Lip}_\sigma$ and $\varkappa$.
For all sufficiently large integers $N$,
%
\begin{equation}
\frac{1}{C_1N^2k^{6\gamma}} \le
\exp\biggl[- \frac{128 t(\mathrm{Lip}_\sigma\vee1)^4}{
\varkappa} k^3\bigl(1+o(1)\bigr)\biggr],
\end{equation}
and the proposition follows upon setting
$\lambda:= \tilde{D} k^{3\gamma}/4$.
\end{pf}

\section{Localization} \label{seclocalization}

The next step in the proof of Theorem~\ref{thmain}
requires us to show that if $x$ and $x'$ are $O(1)$ apart, then $u_t(x)$
and $u_t(x')$ are approximately independent. We show this
by coupling $u_t(x)$ first to the solution of a localized version---see
(\ref{eqheatlocal}) below---of the stochastic heat
equation (\ref{heat1}). And then a second coupling to a suitably-chosen
Picard-iteration approximation of the mentioned localized version.

Consider the following parametric family of random evolution
equations (indexed by the parameter $\beta>0$):
%
\begin{equation}\label{eqheatlocal}\quad
U^{(\beta)}_t(x)
= 1 + \int_{(0,t)\times[ x- \sqrt{\beta t},x+ \sqrt{\beta
t}]}
p_{t-s}(y-x)\sigma\bigl( U^{(\beta)}_s(y)\bigr)
W(\d s \,\d y)
\end{equation}
for all $x\in\R$ and $t\geq0$. 
%
\begin{lemma}\label{lemmomentsUbeta}
Choose and fix $\beta>0$.
Then, (\ref{eqheatlocal}) has an almost surely unique solution
$U^{(\beta)}$ such that for all $T>0$ and $k\ge1$,
%
\begin{equation}
\sup_{\beta>0}\sup_{t\in[0,T]}\sup_{x\in\R}
\E\bigl( \bigl| U^{(\beta)}_t(x)\bigr|^k\bigr) \le
C^k{\mathrm e}^{ak^3},
\end{equation}
where $a$ and $C$ are defined in Lemma~\ref{lemmomentsU}.
\end{lemma}
\begin{pf}
A fixed-point argument shows that there
exists a unique, up to modification, solution to (\ref{eqheatlocal})
subject to the condition that for all $T>0$,
%
\begin{equation}
\sup_{t\in[0,T]}\sup_{x\in\R}\E\bigl(\bigl|U^{(\beta)}_t(x)\bigr|^k
\bigr)<\infty\qquad
\mbox{for all $k\ge1$}.
\end{equation}
See Foondun and Khoshnevisan
\cite{FK} for more details on the ideas of the proof; and the moment estimate
follows as in the proof of Lemma~\ref{lemmomentsU}. We omit
the numerous remaining details.
\end{pf}
%
\begin{lemma}\label{lemuUbeta}
For every $T>0$ there exists a finite and positive constant
$C:=C(\varkappa)$
such that for sufficiently large $\beta>0$ and $k\ge1$,
%
\begin{equation}
\sup_{t\in[0,T]}
\sup_{x\in\R}\E\bigl( \bigl| u_t(x)-U^{(\beta)}_t(x)
\bigr|^k\bigr)
\le C^k k^{k/2}{\mathrm e}^{Fk(k^2 - \beta)},
\end{equation}
where $F\in(0,\infty)$ depends on $(T,\varkappa)$ but not on
$(k
,\beta)$.
\end{lemma}
\begin{pf}
For all $x\in\R$ and $t>0$, define
%
\begin{equation}
V_t(x):= 1 + \int_{(0,t)\times\R} p_{t-s}(y-x)
\sigma\bigl(U^{(\beta)}_s(y)\bigr)
W(\d s \,\d y).\vadjust{\goodbreak}
\end{equation}
Then,
%
\begin{eqnarray}
&&\bigl\| V_t(x)-U^{(\beta)}_t(x)\bigr\|_k\nonumber\\
&&\qquad= \biggl\| \int_{(0,t)\times
\{ y\in\R\dvtx |y-x|>\sqrt{\beta t}\}}
p_{t-s}(y-x)\sigma\bigl(U^{(\beta)}_s(y)\bigr)
W(\d s \,\d y)\biggr\|_k\\
&&\qquad\le2\sqrt k\biggl\| \int_0^t\d s
\int_{|y-x|\ge\sqrt{\beta t}}
\d y\, p^2_{t-s}(y-x) \sigma^2\bigl( U^{(\beta)}_s(y)\bigr)
\biggr\|_{k/2}^{1/2}.\nonumber
\end{eqnarray}
The preceding hinges on an application of
Burkholder's inequality, using the Carlen--Kree bound~\cite{CK}
on Davis's optimal constant~\cite{Davis} in the Burkholder--Davis--Gundy
inequality~\cite{Burkholder,BDG,BG};
see Foondun and Khoshnevisan~\cite{FK} for the details of
the derivation of such an estimate.
Minkowski's inequality tells us then that
the preceding quantity is at most
%
\begin{eqnarray}\label{eqsqsq}\quad
&&2\sqrt{k\int_0^t\d s\int_{|y-x|\ge\sqrt{\beta t}} \d y\,
p_{t-s}^2(y-x) \bigl\| \sigma^2\bigl(
U^{(\beta)}_s(y)\bigr)\bigr\|_{k/2}}\nonumber\\[-8pt]\\[-8pt]
&&\qquad\le\mathrm{const}\cdot
\sqrt{k\int_0^t\d s\int_{|y-x|\ge
\sqrt{\beta t}}\d y\, p_{t-s}^2(y-x)\bigl(1+\bigl\|
U^{(\beta)}_s(y)\bigr\|_{k}^2\bigr).
}\nonumber
\end{eqnarray}
Equation (\ref{eqsqsq}) holds because the Lipschitz continuity of
the function $\sigma$ ensures that it has at-most-linear growth:
$|\sigma(x)|\le\mathrm{const}\cdot(1+|x|)$ for all $x\in\R$.
The inequality in Lemma~\ref{lemmomentsUbeta}
implies that, uniformly over all $t\in[0,T]$ and $x\in\R$,
%
\begin{eqnarray}
\bigl\| V_t(x)-U^{(\beta)}_t(x)\bigr\|_k &\le&
\mathrm{const}\cdot\sqrt{kC^2{\mathrm e}^{2ak^2}
\int_0^t\d r\int_{|z|\ge\sqrt{\beta t}}
\d z\, p^2_r(z)}\nonumber\\[-8pt]\\[-8pt]
&\le&\mathrm{const}\cdot\frac{k^{1/2}{\mathrm e}^{ak^2}}{\sqrt
{\varkappa}}
\sqrt{
\int_0^1\frac{\d s}{\sqrt s}\int_{|w|\ge
\sqrt{2\beta}} \d w\, p_s(w)},\nonumber
\end{eqnarray}
where we have used (\ref{p}). Now a standard Gaussian tail estimate yields
%
\begin{equation}\label{eqgauss}
\int_{|w|\ge\sqrt{2\beta}}p_s(w) \,\d w \le2 {\mathrm e}^{-\beta
/s\varkappa},
\end{equation}
and the latter quantity
is at most $2\exp(-\beta/\varkappa)$ whenever $s\in(0,1]$. Therefore,
on one hand,
%
\begin{equation}\label{eqV-U}
\sup_{x\in\R}
\bigl\| V_t(x)-U^{(\beta)}_t(x)\bigr\|_k
\le\mathrm{const}\cdot
\frac{k^{1/2}{\mathrm e}^{ak^2}}{\sqrt{\varkappa}} {\mathrm e}^{-\beta
/2\varkappa}.
\end{equation}
On the other hand,
\[
u_t(x)-V_t(x) = \int_{(0,t)\times\R} p_{t-s}(y-x)
\bigl[\sigma(u_s(y))-
\sigma\bigl(U^{(\beta)}_s(y)\bigr)\bigr] W(\d s \,\d y),
\]
whence
%
\begin{eqnarray}
&&\| u_t(x)-V_t(x)\|_k\nonumber\\
&&\qquad\le2\sqrt k\biggl\| \int_0^t\d s\int_{-\infty}^\infty\d y\,
p^2_{t-s}(y-x)\bigl[\sigma(u_s(y))-
\sigma\bigl(U^{(\beta)}_s(y)\bigr)\bigr]^2\biggr\|_{k/2}^{1/2}\nonumber\\[-8pt]\\[-8pt]
&&\qquad\le2\sqrt k \operatorname{Lip}_\sigma\biggl\|
\int_0^t\d s\int_{-\infty}^\infty\d y\,
p^2_{t-s}(y-x)\bigl[u_s(y)-U^{(\beta)}_s(y)\bigr]^2
\biggr\|_{k/2}^{1/2}\nonumber\\
&&\qquad\le2\sqrt k \operatorname{Lip}_\sigma\cdot\sqrt{%
\int_0^t\d s\int_{-\infty}^\infty\d y\,
p^2_{t-s}(y-x)\bigl\|
u_s(y)-U^{(\beta)}_s(y)\bigr\|_k^2.}\nonumber
\end{eqnarray}
Consequently, (\ref{eqV-U}) implies that
%
\begin{eqnarray}
&&\bigl\| u_t(x)-U^{(\beta)}_t(x)\bigr\|_k\nonumber\\
&&\qquad\le2\sqrt k \operatorname{Lip}_\sigma\cdot
\sqrt{\int_0^t\d s\int_{-\infty}^\infty\d y\,
p^2_{t-s}(y-x)\bigl\| u_s(y)-U^{(\beta)}_s(y)\bigr\|_k^2}\\
&&\qquad\quad{} +\mathrm{const}\cdot\frac{k^{1/2}
{\mathrm e}^{ak^2}}{\sqrt{\varkappa}} {\mathrm e}^{-\beta/(2\varkappa)}.\nonumber
\end{eqnarray}

Let us introduce a parameter $\delta>0$ and define the seminorms
%
\begin{equation}
\mathcal{N}_{k,\delta}(Z) := \sup_{s\ge0}
\sup_{y\in\R}[
{\mathrm e}^{-\delta s}\| Z_s(y)\|_k]
\end{equation}
for every space--time random field $Z:=\{Z_s(y)\}_{s>0,y\in\R}$.
Then, we have
%
\begin{eqnarray}
&&\mathcal{N}_{k, \delta}\bigl( u-U^{(\beta)}\bigr)\nonumber\\
&&\qquad\le2\sqrt{k} \operatorname{Lip}_\sigma
\mathcal{N}_{k,\delta}
\bigl( u-U^{(\beta)}\bigr)\cdot\sqrt{\int_0^\infty
{\mathrm e}^{-2\delta r} \|p_r\|_{L^2(\R)}^2 \,\d r}\\
&&\qquad\quad{} +\mathrm{const}\cdot
\frac{k^{1/2}{\mathrm e}^{ak^2-\beta/(2\varkappa)}}{\sqrt{\varkappa
}}.\nonumber
\end{eqnarray}
Thanks to (\ref{eqL2pt1}),
if $\delta:= Dk^2$ for some sufficiently large constant $D$, then the
square root
is at most $[4\sqrt{k}(\mathrm{Lip}_\sigma\vee1)]^{-1}$, whence it
follows that (for that fixed choice of $\delta$)
%
\begin{equation}
\mathcal{N}_{k,\delta}\bigl( u-U^{(\beta)}\bigr)
\le\mathrm{const}\cdot
\frac{k^{1/2}{\mathrm e}^{ak^2-\mathrm{const}
\cdot(\beta/\varkappa)}}{\sqrt{\varkappa}}.
\end{equation}
The lemma follows from this.
\end{pf}

Now let us define $U^{(\beta,n)}_t(x)$ to be the $n$th Picard-iteration
approximation
to $U^{(\beta)}_t(x)$. That is, $U^{(\beta,0)}_t(x) := 1$, and for all
$\ell\ge0$,
%
\begin{eqnarray}\label{eqUell}\quad
&&U^{(\beta,\ell+1)}_t(x)\nonumber\\[-8pt]\\[-8pt]
&&\qquad:= 1 + \int_{(0,t)\times[ x-\sqrt{\beta t},
x+\sqrt{\beta t}]} p_{t-s}(y-x) \sigma\bigl( U^{(\beta,\ell
)}_s(y)\bigr)
W(\d s \,\d y).\nonumber
\end{eqnarray}

\begin{lemma}\label{lemiidapprox}
There exist positive and finite constants $C_*$ and
$G$---depend\-ing on $(t,\varkappa)$---such that uniformly for all
$k\in[2,\infty)$ and $\beta>{\mathrm e}$,
%
\begin{equation}
\sup_{x\in\R}
\E\bigl( \bigl| u_t(x) - U^{(\beta,[\log\beta]+1)}_t(x)
\bigr|^k \bigr)
\le\frac{C_*^k k^{k/2} {\mathrm e}^{G k^3}}{\beta^k}.
\end{equation}
\end{lemma}
\begin{pf}
The method of Foondun and Khoshnevisan~\cite{FK} shows
that if $\delta:=D'k^2$ for a sufficiently-large $D'$, then
%
\begin{equation}\label{eqPicardapprox}
\mathcal{N}_{k,\delta}\bigl( U^{(\beta)}-U^{(\beta,n)}\bigr) \le
\mathrm{const}\cdot{\mathrm e}^{-n}\qquad
\mbox{for all $n\ge0$ and $k\in[2,\infty)$}.\hspace*{-35pt}
\end{equation}
%
To elaborate, we follow the arguments in~\cite{FK} of the
proof of Theorem 2.1 leading up to equation (4.6) but with $v_n$ there
replaced by $U^{(\beta,n)}$ here. We then obtain
%
\begin{equation}
\bigl\| U^{(\beta,n+1)}-U^{(\beta,n)} \bigr\|_{k,\theta}
\le\mathrm{const}\cdot\sqrt{k\Upsilon
\biggl(\frac{2\theta}{k}\biggr)}
\bigl\| U^{(\beta,n)}-U^{(\beta,n-1)} \bigr\|_{k,\theta},\hspace*{-35pt}
\end{equation}
where
%
\begin{equation}
\|f\|_{k,\theta}:=\Bigl\{
\sup_{t\ge0}\sup_{x\in\mathbf{R}}
{\mathrm e}^{-\theta t} \E( |f(t,x) |^k )\Bigr\}^{1/k}
\end{equation}
and
%
\begin{equation}
\Upsilon(\theta):=\frac{1}{2\pi}
\int_{-\infty}^{\infty} \frac{\d\xi}{\theta+\xi^2}.
\end{equation}
A quick computation reveals that by choosing
$\theta:= D''k^3$, for a large enough constant $D''>0$, we obtain
%
\begin{equation}
\bigl\| U^{(\beta,n+1)}-U^{(\beta,n)} \bigr\|_{k,\theta}
\le{\mathrm e}^{-1} \bigl\| U^{(\beta,n)}-U^{(\beta,n-1)} \bigr\|
_{k,\theta}.
\end{equation}
We get (\ref{eqPicardapprox}) from this.

Next we set $n:=[\log\beta]+1$ and apply the preceding together with Lem\-ma~\ref{lemuUbeta} to finish the proof.
\end{pf}
%
\begin{lemma}\label{lemiid}
Choose and fix $\beta,t>0$ and $n\ge0$. Also fix $x_1,x_2,\ldots\in
\R$
such that $|x_i-x_j| \ge2n\sqrt{\beta t}$ whenever $i\neq j$.
Then
$\{U^{(\beta,n)}_t(x_j)\}_{j\in\Z}$ is a collection of
i.i.d. random variables.
\end{lemma}
\begin{pf}
The proof uses induction on the variable $n$,
and proceeds by establishing a little more. We will use the
$\sigma$-algebras $\mathcal{P}(A)$\vadjust{\goodbreak}
as defined in Appendix~\ref{appSI}, where $A\subset\R$
ranges over all Lebesgue-measurable sets of finite Lebesgue measure
Proposition~\ref{prmeasWalshint}.

Since $U^{(\beta,0)}_t(x) \equiv1$, the statement of the lemma holds
tautologically
for $n=0$. In order to understand the following argument better let us
concentrate
on the case $n=1$. In that case,
%
\begin{equation}\quad
U^{(\beta,1)}_t(x) = 1 + \sigma(1) \cdot\int_{(0,t)\times[
x-\sqrt{\beta t},x+\sqrt{\beta t}]} p_{t-s}(y-x)
W(\d s \,\d y).
\end{equation}
In particular, if we define the process
$\{\tilde{U}^{(\beta,1)}_s(t,x)\}_{0\leq s \leq t}$ by
%
\begin{equation}\qquad
\tilde{U}^{(\beta,1)}_s(t,x) = 1 + \sigma(1) \cdot\int
_{(0,s)\times
[x-\sqrt{\beta t},x+\sqrt{\beta t}]} p_{t-r}(y-x)
W(\d r \,\d y),
\end{equation}
it follows from Proposition~\ref{prmeasWalshint}
that $U^{(\beta,1)}_{\bullet}(t,x)\in
\mathcal{P}([x-\sqrt{\beta t},x+\sqrt{\beta t}])$. Hence,
Corollary~\ref{comeasWalshint} shows that
$\{\tilde{U}^{(\beta,1)}_{\bullet}(t,x_j)\}_{j\in\Z}$ are independent
processes. Taking $s=t$ shows the independence part of the lemma for $n=1$.
It is not hard to see that the law of $U^{(\beta,1)}_t(x)$ is
independent of $x$,
namely a Gaussian distribution with mean and variance
parameters that are independent of $x$.
This concludes the proof for $n=1$.

Now, let us define processes $\{\tilde{U}^{(\beta,n)}_s(t,x)\}_{0
\leq
s \leq t}$ by
%
\begin{eqnarray}
&&\tilde{U}^{(\beta,n)}_s(t,x) \nonumber\\
&&\qquad= 1 + \sigma(1) \cdot\int_{(0,s)\times[
x-\sqrt{\beta t},x+\sqrt{\beta t}]} p_{t-r}(y-x)
\sigma\bigl(\tilde{U}^{(\beta,n-1)}_r(r,y)\bigr)\\
&&\hspace*{170pt}{}\times W(\d r \,\d y).\nonumber
\end{eqnarray}
If we proved that
$\tilde{U}^{(\beta,n)}_{\bullet}(t,x)\in\mathcal{P}([x-n\sqrt
{\beta
t},
x+n\sqrt{\beta t}])$ for all $x\in\R$ and $t>0$,
it then would follow from (\ref{eqUell})
and Proposition~\ref{prmeasWalshint} that
$\tilde{U}^{(\beta,n+1)}_{\bullet}(t,x)
\in\mathcal{P}([x-(n+1)\sqrt{\beta t},x+(n+1)\sqrt{\beta t}])$
for all $x\in\R$ and $t>0$. Since this fact is true for $n=1$,
then we have proved that
$\tilde{U}^{(\beta,n)}_{\bullet}(t,x)\in\mathcal{P}([x-n\sqrt
{\beta
t},
x+n\sqrt{\beta t}])$ for all $x\in\R$, $t>0$ and $n\in\N$.
We then use this result, together with Corollary~\ref{comeasWalshint}
to deduce that
$\{\tilde{U}^{(\beta,n)}_{\bullet}(t,x_j)\}_{j\in\Z}$ are independent
processes. The independence part of the lemma follows upon taking
$s:=t$ in this discussion.

Since the law of the noise $W$ is invariant by translation,
it is not hard to prove by induction that the law of $U^{(\beta,n)}_t(x)$
is indeed independent of $x$.
(Notice that this is always true when the initial condition is constant
\cite{Dalang99}, Lemma 18.) This concludes the (inductive) proof of
the lemma.
\end{pf}

\section{\texorpdfstring{Proof of Theorem \protect\ref{thmain}}{Proof of Theorem 1.1}} \label{secproofmain}

We are ready to combine our efforts thus far in order to verify
Theorem~\ref{thmain}.
\begin{pf*}{Proof of Theorem~\ref{thmain}}
Parts (1) and (2) of Theorem~\ref{thmain} are proved similarly. Therefore,
we present the details of the second part.
For the proof of the first part, we can take $\gamma:=1/6$
in the following argument. We remind that the processes $U^{\beta}$
and $U^{(\beta,n)}$ are defined, respectively, in (\ref{eqheatlocal}) and
(\ref{eqUell}).

For all $x_1,\ldots,x_N\in\R$,
%
\begin{eqnarray}
\P\Bigl\{ {\max_{1\le j\le N}} |u_t(x_j)|<\lambda\Bigr\} &\le&
\P\Bigl\{ \max_{1\le j\le N}\bigl|U_t^{(\beta, \log\beta)}(x_j)\bigr|
<2\lambda\Bigr\} \nonumber\\[-8pt]\\[-8pt]
&&{} + \P\Bigl\{ \max_{1\le j\le N}\bigl|
U_t^{(\beta,\log\beta)}(x_j)-u_t(x_j)\bigr| >\lambda
\Bigr\}.\nonumber
\end{eqnarray}
(To be very precise, we need to write $(\beta,[\log\beta]+1)$ in
place of
$(\beta,\log\beta)$.)
The whole program of Section~\ref{sectailestimates}, that led to
Proposition~\ref{prsigmadecay} can be carried out for
$U^{(\beta)}$ instead of $u$. Only minor changes
(typically in Proposition~\ref{prlowermoments}) are needed.
Then, (\ref{eqPicardapprox}) shows that the same moment estimates
are valid for $U^{(\beta,n)}$ as for $U^{(\beta)}$.

We can follow along the proof of Proposition~\ref{prsigmadecay}
and prove similarly
the existence of constants $c_1,c_2>0$---independent
of $\beta$ for all sufficiently large values of $\beta$---so that
for all $x\in\R$ and $\lambda\ge1$,
%
\begin{equation}
\label{eqUbetalbd}
\P\bigl\{\bigl\vert U_t^{(\beta,\log\beta)}(x)\bigr\vert\ge\lambda
\bigr\}
\ge c_1{\mathrm e}^{-c_2\lambda^{1/\gamma}}.
\end{equation}

Suppose, in addition, that $|x_i-x_j|\ge2\sqrt{\beta t}\log\beta$ whenever
$i\neq j$. Then, Lemmas~\ref{lemiidapprox} and~\ref{lemiid}
together imply the following:
%
\begin{eqnarray} \label{equbd}
&& \P\Bigl\{ {\max_{1\le j\le N}}
|u_t(x_j)|<\lambda\Bigr\}\nonumber\\[-8pt]\\[-8pt]
&&\qquad\le\bigl(1-c_1{\mathrm e}^{-c_2\cdot(2\lambda)^{1/\gamma
}}\bigr)^N
+N C_*^k k^{k/2}{\mathrm e}^{Gk^3}
\beta^{-k}\lambda^{-k}.\nonumber
\end{eqnarray}
The constants $C_{*}$ and $G$ may differ from the ones in
Lemma~\ref{lemiidapprox}. We now select the various parameters judiciously:
choose $\lambda:= k$, $N := \lceil k\exp(c_2\cdot(2k)^{1/\gamma
})\rceil$ and
$\beta:= \exp(\rho k^{(1-\gamma)/\gamma})$
for a large-enough positive constant $\rho>2\cdot3^{1/\gamma}c_2$.
In this way, (\ref{equbd}) simplifies:
for all sufficiently large integers $k$,
%
\begin{eqnarray}
&& \P\Bigl\{ {\max_{1\le j\le N}}
|u_t(x_j)|<k\Bigr\}\nonumber\\
&&\qquad\le{\mathrm e}^{-c_1k}+\exp\biggl[c_2\cdot(2k)^{1/\gamma}+
\log k + k\log C_* -\frac{k\log k}{2} +
Gk^3-\frac{\rho k^{1/\gamma}}{2} \biggr]\hspace*{-10pt}\\
&&\qquad\le2{\mathrm e}^{-c_1k}.\nonumber
\end{eqnarray}

Now we choose the $x_i$'s as follows: set $x_0:=0$,
and define iteratively
%
\begin{eqnarray}\label{eqxi}
x_{i+1} :\!&=& x_i + 2\sqrt{\beta t}([\log\beta]+1)
\nonumber\\[-8pt]\\[-8pt]
&=& 2(i+1)\sqrt{\beta t}
( [\log\beta]+1)\qquad
\mbox{for all $i\ge0$}.\nonumber
\end{eqnarray}
The preceding implies that for all sufficiently large $k$,
%
\begin{equation}
\P\Bigl\{ {\sup_{x\in[0,2(N+1)\sqrt{\beta t}
( [\log\beta]+1)]}} |u_t(x)| <k\Bigr\}
\le2{\mathrm e}^{-c_1k},
\end{equation}
whence
%
\begin{equation}
\P\Bigl\{ {\sup_{|x|\le2(N+1)\sqrt{\beta t}
( [\log\beta]+1)}} |u_t(x)| <k\Bigr\}
\le2{\mathrm e}^{-c_1k}
\end{equation}
by symmetry.

It is not hard to verify that, as $k\to\infty$,
%
\begin{equation}
2(N+1) \sqrt{\beta t}([\log\beta]+1) =
O( {\mathrm e}^{\rho k^{1/\gamma}}).
\end{equation}
Consequently,
the Borel--Cantelli lemma, used in conjunction with a standard monotonicity
argument, implies that $u_t^*(R)\ge\mathrm{const}\cdot
(\log(R)/c_2)^{\gamma}$
a.s. for all sufficiently large values of $R$, where $u_t^*(R)$ is
defined in (\ref{equstar}). By Proposition~\ref{prsigmadecay},
$c_2 = \mathrm{const} \cdot\varkappa^{1/12\gamma}$.
Therefore, the theorem follows.
\end{pf*}

\section{\texorpdfstring{Proof of Theorem \protect\ref{thboundedsigma}}{Proof of Theorem 1.2}}
\label{secboundedsigma}

Next we prove Theorem~\ref{thboundedsigma}.
In order to obtain the upper bound, the proof requires an
estimate of spatial continuity of
$x\mapsto u_t(x)$. However, matters are somewhat complicated
by the fact that we need a modulus of continuity estimate
that holds simultaneously for every
$x\in[-R,R]$, uniformly for all large values of $R$. This will be
overcome in a few steps. The first is a standard moment bound
for the increments of the solution;
however, we need to pay close attention to the constants in the
estimate.\looseness=-1

\begin{lemma}\label{lemMC1}
Choose and fix some $t>0$, and suppose $\sigma$
is uniformly bounded.
Then there exists a finite and positive constant $A$ such that
for all real numbers $k\ge2$,
%
\begin{equation}
\sup_{-\infty<x\neq x'<\infty}
\frac{\E(|u_t(x)-u_t(x')|^{2k})}{|x-x'|^k}\le
\biggl(\frac{Ak}{\varkappa}\biggr)^k.
\end{equation}
\end{lemma}
\begin{pf}
Throughout, let $S_0:={\sup_{x\in\R} }|\sigma(x)|$.

If $x,x'\in[-R,R]$ and $t>0$ are held fixed, then
we can write $u_t(x) - u_t(x') = N_t$,
where $\{N_\tau\}_{\tau\in(0,t)}$ is the continuous
mean-one martingale described by
%
\begin{equation}
N_\tau:=\int_{(0,\tau)\times\R}
[ p_{t-s}(y-x)-p_{t-s}(y-x')]\sigma(u_s(y)) W(\d s \,\d y)
\end{equation}
for $\tau\in(0,t)$. The quadratic variation of $\{N_\tau\}_{\tau
\in(0,t)}$
is estimated as follows:
%
\begin{eqnarray}
\langle N\rangle_\tau&\le& S_0^2\cdot
\int_0^\tau\d s\int_{-\infty}^\infty\d y\,
[ p_s(y-x)-p_s(y-x')]^2\nonumber\\[-8pt]\\[-8pt]
&\le& {\mathrm e}^\tau S_0^2\cdot\int_0^\infty{\mathrm e}^{-s} \,\d s\int
_{-\infty}^\infty
\d y\,
[ p_s(y-x)-p_s(y-x')]^2.\nonumber
\end{eqnarray}
For every $s>0$ fixed, we can compute the $\d y$-integral
using Plancherel's theorem, and obtain
$\pi^{-1}\int_{-\infty}^\infty
(1-\cos(\xi|x-x'|)) \exp(-\varkappa s\xi^2) \,\d\xi$.
Therefore, there exists a finite and positive constant $a$ such that
%
\begin{equation} \label{eqnquadvar}
\langle N\rangle_\tau\le\frac{{\mathrm e}^\tau S_0^2}{\pi}\cdot
\int_{-\infty}^\infty
\frac{1-\cos(\xi|x-x'|)}{1+\varkappa\xi^2} \,\d\xi\\
\le\frac{a}{\varkappa} |x-x'|,
\end{equation}
uniformly for all $\tau\in(0,t)$;
we emphasize that $a$ depends only on $S_0$ and~$t$.
The Carlen--Kree estimate~\cite{CK} for the Davis~\cite{Davis}
optimal constant in the Burkholder--Davis--Gundy
inequality~\cite{Burkholder,BDG,BG} implies the lemma.
\end{pf}

The second estimate turns the preceding moment bounds into
an maximal exponential estimate. We use a standard chaining argument
to do this. However, once again we have to pay close attention to
the parameter dependencies in the implied constants.
%
\begin{lemma}\label{lemMC2}
Choose and fix $t>0$, and suppose $\sigma$ is
uniformly bounded. Then there exist a constant
$C\in(0,\infty)$ such that
%
\begin{equation}
\E\biggl[\mathop{\sup_{x,x'\in I:}}_{|x-x'|\le\delta}\exp
\biggl(\frac{
\varkappa|u_t(x)-u_t(x')|^2}{
C\delta} \biggr) \biggr] \le\frac{2}{\delta},
\end{equation}
uniformly for every $\delta\in(0,1]$
and every interval $I\subset[0,\infty)$ of length at most one.
\end{lemma}
\begin{pf}
Recall~\cite{Minicourse}, (39), page 11, the Kolmogorov
continuity in the following quantitative form:
suppose there exist $\nu>\gamma>1$ for which
a stochastic process
$\{\xi(x)\}_{x\in\R}$ satisfies the following:
%
\begin{equation}
\E\bigl(|\xi(x)-\xi(x')|^\nu\bigr) \le C|x-x'|^\gamma;
\end{equation}
we assume that the preceding holds
for all $x,x'\in\R$, and $C\in(0,\infty)$ is
independent of $x$ and $x'$. Then, for every integer $m\ge0$,
%
\begin{equation}
\E\Bigl({\mathop{\sup_{x,x'\in I:}}_{|x-x'|\le2^{-m}}}
|\xi(x)-\xi(x')|^\nu\Bigr) \le\biggl(
\frac{2^{(2-\gamma+\nu)/\nu} C^{1/\nu}}{1-2^{-(\gamma-1)/\nu}}
\biggr)^\nu\cdot2^{-m(\gamma-1)}.
\end{equation}
(Reference~\cite{Minicourse}, (39), page 11, claims this with
$2^{-m(\gamma-1)}$
replaced with $2^{-m\gamma}$
on the right-hand side. But this is a typographical error; compare
with~\cite{Minicourse}, (38), page~11.)

If $\delta\in(0,1]$, then we can find an integer $m\ge0$
such that $2^{-m-1}\le\delta\le2^{-m}$, whence
it follows that
%
\begin{eqnarray}
\E\Bigl({\mathop{\sup_{x,x'\in I:}}_{|x-x'|\le\delta}}
|\xi(x)-\xi(x')|^\nu\Bigr) &\le&
\biggl(\frac{2^{(2-\gamma+\nu)/\nu} C^{1/\nu}}{1-2^{-(\gamma-1)/\nu}}
\biggr)^\nu\cdot2^{-m(\gamma-1)}\nonumber\\[-8pt]\\[-8pt]
&\le&
\biggl(\frac{2^{(2-\gamma+\nu)/\nu} C^{1/\nu}}{1-2^{-(\gamma-1)/\nu}}
\biggr)^\nu\cdot(2\delta)^{\gamma-1}.\nonumber
\end{eqnarray}
We apply the preceding with $\xi(x):=u_t(x)$,
$\gamma:=\nu/2:=k$ and $C:=(Ak/\varkappa)^{k}$, where $A$ is the constant
of Lemma~\ref{lemMC1}. It follows that there exists a positive
and finite constant $A_*$ such that for all
intervals $I$ of length at most one, all integers $k\ge2$,
and every $\delta\in(0,1]$,
%
\begin{equation}
\E\Bigl({\mathop{\sup_{x,x'\in I:}}_{|x-x'|\le\delta}}
|u_t(x)-u_t(x')|^{2k}\Bigr) \le\biggl(\frac{A_*k}{\varkappa}
\biggr)^k \delta^{k-1}.
\end{equation}
Stirling's\vspace*{1pt} formula tells us that
there exists a finite constant $B_* > 1$ such that
$(A_*k)^k\le B_*^k k!$ for all integers $k\ge1$. Therefore,
for all $\alpha,\delta>0$,
%
\begin{equation} \label{eq610}
\E\biggl[\mathop{\sup_{x,x'\in I:}}_{|x-x'|\le\delta}\exp
\biggl(\frac{
\alpha|u_t(x)-u_t(x')|^2}{
\delta} \biggr) \biggr] \le
\frac1\delta\sum_{k=0}^\infty\biggl(\frac{\zeta B_*}{\varkappa
}\biggr)^k.
\end{equation}
And this is at most two if $\alpha:=\varkappa/(2B_*)$.
The result follows.
\end{pf}

Next we obtain another moments bound, this time for the solution
rather than its increments.
%
\begin{lemma}\label{lemMC3}
Choose and fix $t>0$, and suppose $\sigma$ is uniformly bounded.
Then for all integers $k\ge1$,
%
\begin{equation}
\sup_{x\in\R} \E(|u_t(x)|^{2k})
\le\bigl(2\sqrt{2}+o(1)\bigr) (\tilde\mu_t k)^{k}\qquad
(\mbox{as $k\to\infty$})
\end{equation}
where the ``$o(1)$'' term depends only on $k$, and
%
\begin{equation}
\tilde{\mu}_t := \frac{8}{{\mathrm e}}\cdot S_0^2\sqrt{\frac{t}{\pi
\varkappa}}.
\end{equation}
\end{lemma}
\begin{pf}
Let us choose and fix a $t>0$.
Define $S_0:={\sup_{x\in\R}}|\sigma(x)|$, and recall the
martingale $\{M_\tau\}_{\tau\in(0,t)}$ from
(\ref{eqMtau}). It\^o's formula (\ref{eqIto}) tells us that
a.s., for all $\tau\in(0,t)$,
%
\begin{equation}\qquad
M_\tau^{2k} \le1+2k\int_0^\tau M_s^{2k-1} \,\d M_s+
\pmatrix{2k\cr2} S_0^2\int_0^\tau\frac{M_s^{2(k-1)}}{
(4\pi\varkappa(t-s))^{1/2}} \,\d s.
\end{equation}
[Compare with (\ref{eqItoLB}).] We can take expectations,
iterate the preceding and argue as we did in the proof of Proposition
\ref{prlowermoments}. To summarize the end result, let us define
%
\begin{equation}\quad
\eta_t(x):=1+S_0\cdot\int_{(0,\tau)\times
\mathbf{R}} p_{t-s}(y-x) W(\d s \,\d y)\qquad (0\le\tau\le t)
\end{equation}
and
%
\begin{equation}
\zeta_t(x):=S_0\cdot\int_{(0,\tau)\times
\mathbf{R}} p_{t-s}(y-x) W(\d s \,\d y)\qquad (0\le\tau\le t).
\end{equation}
Then we have
$\E[M_t^{2k} ]\le\E[\eta_t(x)^{2k} ]\le
2^{2k} (1+ \E[\zeta_t(x)^{2k} ] )$, and
similar computations as those
in the proof of Proposition~\ref{prlowermoments} yield the
lemma.\vspace*{-2pt}
\end{pf}

Next we turn the preceding moment bound into a sharp Gaussian
tail-probability estimate.\vspace*{-2pt}
%
\begin{lemma}\label{lemMC4}
Choose and fix a $t>0$, and suppose that $\sigma$ is uniformly
bounded. Then there exist finite constants $C>c>0$ such that
simultaneously for all $\lambda>1$ and $x\in\R$,
%
\begin{equation}
c\exp\bigl(-C
\sqrt{\varkappa} \lambda^2\bigr)
\le\P\{|u_t(x)|>\lambda\} \le C \exp\bigl(-
c \sqrt{\varkappa} \lambda^2\bigr).\vspace*{-2pt}
\end{equation}
\end{lemma}
\begin{pf}
The lower bound is proved by an appeal to
the Paley--Zygmund inequality, in the very
same manner that Proposition~\ref{prlowertail}
was established. However, we apply the improved
inequality in Lemma~\ref{lemMC3}
(in place of the result of Lemma~\ref{lemmomentsU}).
As regards the upper bound, note that Lemma~\ref{lemMC3}
implies that there exists a positive and finite constant
$\tilde{A}$ such that for all integers $m\ge0$,
$\sup_{x\in\R}\E(|u_t(x)|^{2m})
\le(\tilde{A}/\sqrt{\varkappa})^m m!$,
thanks to the Stirling formula. Thus,
%
\begin{equation}
\sup_{x\in\R}\E\exp(\alpha|u_t(x)|^2)
\le\sum_{m=0}^\infty\biggl(\frac{\alpha\tilde{A}}{%
\sqrt{\varkappa}}\biggr)^m
=\frac{1}{1-\alpha\tilde{A} \varkappa^{-1/2}}<\infty,
\end{equation}
provided that $\alpha\in(0,\sqrt{\varkappa}/\tilde{A})$.
Notice that this has a different behavior than (\ref{eq610})
in terms of $\kappa$. If we fix such an
$\alpha$, then we obtain from Chebyshev's inequality the bound
$\P\{ u_t(x) >\lambda\} \le
(1-\alpha\tilde{A} \varkappa^{-1/2})^{-1}\cdot\exp(-\alpha\lambda^2)$,
valid simultaneously for all $x\in\R$ and $\lambda>0$.
We write $\alpha:= c\sqrt{\varkappa}$ to finish.\vspace*{-2pt}
\end{pf}

We are finally ready to assemble the preceding estimates in
order to establish Theorem~\ref{thboundedsigma}.\vspace*{-2pt}
\begin{pf*}{Proof of Theorem~\ref{thboundedsigma}}
Consider the proof of Theorem~\ref{thmain}:
if we replace the role of (\ref{eqtails}) by the bounds
in Lemma~\ref{lemMC4} and choose $\lambda:=k$,
$N:=\lceil k \times\exp(c\sqrt{\varkappa} k^2)\rceil$ and
$\beta:= \exp((\mathrm{Lip}_\sigma\vee1)^4 k^2 / \varkappa)$
in the equivalent of (\ref{equbd}) with the appropriate estimates,
then we obtain the almost-sure bound
$\liminf_{R\to\infty} u_t^*(R)/(\log R)^{1/2}>
\mathrm{const}\cdot\varkappa^{-1/4}$, where ``const'' is independent of
$\varkappa$.
It remains to derive a corresponding upper bound
for the $\limsup$.

Suppose\vspace*{1pt} $R\ge1$ is an integer. We partition the interval
$[-R,R]$ using a length-1 mesh with endpoints
$\{x_j\}_{j=0}^{2R}$ via
%
\begin{equation}
x_j := -R+j \qquad\mbox{for $0\le j\le2R$}.
\end{equation}
Then we write
%
\begin{equation} \label{eqprustar}
\P\{ u_t^*(R)> 2\alpha(\log R)^{1/2}\}\le
\mathcal{T}_1+\mathcal{T}_2,\vadjust{\goodbreak}
\end{equation}
where
%
\begin{eqnarray}
\mathcal{T}_1&:=&\P\Bigl\{ \max_{1\le j\le2R} u_t (x_j)
>\alpha(\log R)^{1/2}\Bigr\},\nonumber\\[-9pt]\\[-9pt]
\mathcal{T}_2 &:=& \P\Bigl\{{\max_{1\le j\le2R}
\sup_{x\in(x_j,x_{j+1})}} |u_t(x)-u_t(x_j)|>
\alpha(\log R)^{1/2} \Bigr\}.\nonumber
\end{eqnarray}
By Lemma~\ref{lemMC4},
%
\begin{equation}
\mathcal{T}_1 \le2R\sup_{x\in\R} \P\{ u_t(x)>\alpha(\log
R)^{1/2}\}
\le\frac{\mathrm{const}}{R^{-1+c\sqrt{\varkappa} \alpha^2}}.
\end{equation}
Similarly,
%
\begin{equation}
\mathcal{T}_2 \le2R\sup_{I}\P\Bigl\{
{\sup_{x,x'\in I}}|u_t(x)-u_t(x')|>
\alpha(\log R)^{1/2}\Bigr\},
\end{equation}
where ``$\sup_I$'' designates a supremum over all
intervals $I$ of length one. Chebyshev's inequality and Lemma~\ref{lemMC2}
together imply that
%
\begin{eqnarray}
\mathcal{T}_2 &\le& 2R^{-(\varkappa/C) \alpha^2+1}
\sup_{I}\E\biggl[\sup_{x,x'\in I}\exp\biggl(
\frac{\varkappa}{C}|u_t(x)-u_t(x')|^2\biggr)\biggr]\nonumber\\[-9pt]\\[-9pt]
&\le& \frac{\mathrm{const}}{R^{-1+(\varkappa
\alpha^2)/C}}.\nonumber
\end{eqnarray}
Let $q:=\min(\varkappa/C,c\sqrt{\varkappa})$ to find that
%
\begin{equation}
\sum_{R=1}^\infty
\P\{ u_t^*(R) > 2\alpha(\log R)^{1/2}\}
\le\mathrm{const}\cdot\sum_{R=1}^\infty
R^{-q\alpha^2+1},
\end{equation}
and this is finite provided that $\alpha>(2/q)^{1/2}$.
By the Borel--Cantelli lemma,
%
\begin{equation}
\mathop{\limsup_{R\to\infty:}}_{R\in\Z}
\frac{u_t^*(R)}{(\log R)^{1/2}}\le
\biggl(\frac8q\biggr)^{1/2}
<\infty\qquad\mbox{a.s.}
\end{equation}
Clearly, $(8/q)^{1/2}
\le\mathrm{const} \cdot\varkappa^{-1/4}$ for all $\varkappa\ge
\varkappa_0$,
for a constant depends only on $\varkappa_0$.
And we can remove the restriction ``$R\in\Z$'' in the $\limsup$
by a standard monotonicity argument; namely, we find---by
considering in the following $R-1\le X\le R$---that
%
\begin{equation}\quad
\limsup_{X\to\infty}
\frac{u_t^*(X)}{(\log X)^{1/2}} \le
\mathop{\limsup_{R\to\infty:}}_{R\in\Z}
\frac{u_t^*(R)}{(\log(R-1))^{1/2}}
\le\biggl(\frac8q\biggr)^{1/2} \qquad\mbox{a.s.}
\end{equation}
This proves the theorem.
\end{pf*}

\section{\texorpdfstring{Proof of Theorem \protect\ref{thPAM}}{Proof of Theorem 1.3}}

This section is mainly concerned with the proof of Theorem \ref
{thPAM}. For that purpose, we start with tail-estimates.
%
\begin{lemma} \label{lempatails}
Consider (\ref{heat1}) with $\sigma(x) := cx$, where
$c > 0$ is fixed. Then,
%
\begin{equation}
\log\P\{ |u_t(x)|\ge\lambda\}
\asymp- \sqrt{\varkappa} (\log\lambda)^{3/2}\qquad
\mbox{as $\lambda\to\infty$}.\vadjust{\goodbreak}
\end{equation}
\end{lemma}
\begin{pf}
Corollary~\ref{comomentsU} implies the upper bound
(the boundedness of $\sigma$ is not required in the results of
Section~\ref{secuppertail}).

As for the lower bound, we know from
\cite{BC}, Theorem 2.6,
%
\begin{equation} \label{eqBCmoments}
{\mathrm e}^{k(k^2-1)t /24 \varkappa} \le\E(|u_t(x)|^k)
\le2{\mathrm e}^{k(k^2-1)t / 24\varkappa},
\end{equation}
uniformly for all integers $k\ge2$ and $x \in\R$.
Now we follow the same method as in the proof of Proposition
\ref{prlowertail}, and use the Paley--Zygmund inequality to obtain
%
\begin{eqnarray}
\P\biggl\{|u_t(x)| \ge\frac{1}{2}\|u_t(x)\|_{2k}
\biggr\} &\ge& \frac{[\E(|u_t(x)|^{2k})
]^2}{%
4\E(|u_t(x)|^{4k})} \nonumber\\[-8pt]\\[-8pt]
&\ge& C_1 {\mathrm e}^{-D_1k^3 / \varkappa}\nonumber
\end{eqnarray}
for some nontrivial constants $C_1$ and $D_1$ that do not depend
on $x$ or $k$. We then obtain the following: uniformly for all $x\in\R$
and sufficiently-large integers $k$,
%
\begin{equation}
\P\biggl\{|u_t(x)|\ge\frac{C}{2}
{\mathrm e}^{4Dk^2 / \varkappa}\biggr\}\ge C_1{\mathrm e}^{-D_1k^3 /
\varkappa}.
\end{equation}
Let $\lambda:=(C/2)\exp\{4Dk^2 / \varkappa\}$, and apply
a direct computation to deduce the lower bound.
\end{pf}

We are now ready to prove Theorem~\ref{thPAM}.
Our proof is based on roughly-similar ideas to those used in the course
of the proof
of Theorem~\ref{thmain}. However, at a technical level, they are
slightly different. Let us point out some of the essential differences:
unlike what we did in the proof of Theorem~\ref{thmain},
we now do not choose the values of $N$, $\beta$ and $\lambda$
as functions of $k$, but rather as functions of $R$;
the order of the moments $k$ will be fixed; and we
will not sum on $k$, but rather sum on a discrete sequence of values of
the parameter~$R$. The details follow.
\begin{pf*}{Proof of Theorem~\ref{thPAM}}
First we derive the lower bound by following the same
method that was used in the proof of Theorem~\ref{thmain};
see Section~\ref{secproofmain}. But
we now use Lemma~\ref{lempatails} rather than Corollary~\ref{cobds}.

The results of Section~\ref{seclocalization}
can be modified to apply to the parabolic Anderson model, provided that
we again apply Lemma~\ref{lempatails} in place of Corollary~\ref{cobds}.
In this way we obtain the following, where the $x_i$'s are
defined by (\ref{eqxi})\setcounter{footnote}{3}\footnote{To be very
precise, we once again need to write
$(\beta,[\log\beta]+1)$ in place of
$(\beta,\log\beta)$.}:
consider the event
%
\begin{equation}
\Lambda:= \Bigl\{ {\max_{1\le j\le N}} |u_t(x_j)| <
\Xi\Bigr\}\qquad
\mbox{where } \Xi:=\exp\biggl(
C_1\frac{(\log R)^{2/3}}{\varkappa^{1/3}}
\biggr).
\end{equation}
Then,
%
\begin{eqnarray}\label{eqpam1}
\P(\Lambda)
&\le& \P\Bigl\{ \max_{1\le j\le N}
\bigl|U^{(\beta,\log\beta)}_t(x_j)\bigr| <2\Xi\Bigr\} \nonumber\\[-2pt]
&&{} +\P\bigl\{ \bigl|u_t(x_j) -U_t^{(\beta,\log\beta)}\bigr| >
\Xi\mbox{ for some } 1\le j\le N\bigr\} \\[-2pt]
&\le& \bigl(1-\P\bigl\{ \bigl|U^{(\beta,\log\beta)}_t(x_j)\bigr|
\ge2\Xi\bigr\}\bigr)^N
+ \frac{N\beta^{-k}C_*^kk^{k/2}{\mathrm e}^{Gk^3}}{\Xi}.\nonumber
\end{eqnarray}
Note that we do not yet have a lower bound on
$\P\{|U_t^{(\beta,\log\beta)}(x)|\ge\lambda\}$. However, we have
%
\begin{eqnarray}\label{eqpam2}
&&\P\bigl\{ \bigl|U^{(\beta,\log\beta)}_t(x_j)\bigr| \ge2
\Xi\bigr\}\nonumber \\[-2pt]
&&\qquad\ge\P\{ |u_t(x_j)| \ge3
\Xi\} -\P\bigl\{ \bigl|u_t(x_j)-U^{(\beta,\log\beta)}_t(x_j)\bigr|
\ge\Xi\bigr\} \\[-2pt]
&&\qquad \ge\alpha_1R^{-\alpha_2C_1^{3/2}}-\frac{N\beta^{-k}
C_*^k k^{k/2}{\mathrm e}^{Gk^3}}{\Xi},\nonumber
\end{eqnarray}
valid for some positive constants $\alpha_1$ and $\alpha_2$.
Now let us choose $N:=\lceil R^{a}\rceil$ and $\beta:=R^{1-a}$
for a fixed $a\in(0,1)$. With these values of $N$ and $\beta$
and the lower bound in (\ref{eqpam2}), the upper bound in
(\ref{eqpam1}) becomes
%
\begin{equation}
\P(\Lambda)\le\biggl(1- \alpha_1 R^{-\alpha_2C_1^{3/2}}
+\frac{C_*^kk^{k/2}{\mathrm e}^{Gk^3}}{R^{k(1-a)-a}
\Xi} \biggr)^N
+ \frac{C_*^kk^{k/2}{\mathrm e}^{Gk^3}}{R^{k(1-a)-a}
\Xi}.
\end{equation}
Let us consider $k$ large enough so that
$k(1-a)-a>2$. Notice that $k$ will not depend on $R$; this is in contrast
with what happened in the proof of Theorem~\ref{thmain}.

We can choose the constant $C_1$ to be small enough to satisfy
$ \alpha_2 C_1^{3/2} < a/2$.
Using these, we obtain
%
\begin{equation}
\P\Bigl\{ {\sup_{x\in[0,R]} }|u_t(x)| <
{\mathrm e}^{C_1(\log R)^{2/3}/\varkappa^{1/3}}\Bigr\}
\le\exp(-\alpha_1R^{a/2} )+\frac{\mathrm{const}}{R^2}.
\end{equation}
The Borel--Cantelli lemma yields the lower bound of the theorem.

We can now prove the upper bound. Our derivation is
modeled after the proof of Theorem~\ref{thboundedsigma}.

First, we need a continuity estimate for the solution of
(\ref{heat1}) in the case that $\sigma(x) := cx$. In accord with
(\ref{eqBCmoments}),
%
\begin{eqnarray}\qquad
&&\E\bigl( \vert u_t(x)-u_t(y)\vert^{2k}\bigr)\nonumber\\[-2pt]
&&\qquad\le\bigl( 2\sqrt{2k} \bigr)^{2k}\biggl[
\int_0^{t} \d r\, \| u(r,0)\|_{2k}^2 \int_{\R} \d z\,
| p_{t-r}(x-z)-p_{t-r}(y-z)|^2 \biggr]^{k} \\[-2pt]
&&\qquad\le\bigl( 2\sqrt{2k}\bigr)^{2k}
\biggl[\int_0^{t} \d r\, 2^{1/k}
{\mathrm e}^{8Dk^2/\varkappa} \int_{\R} \d z\, |
p_{t-r}(x-z)-p_{t-r}(y-z)|^2 \biggr]^{k}\nonumber
\end{eqnarray}
for some constant $D$ which depends on $t$. Consequently [see the
derivation of (\ref{eqnquadvar})],
%
\begin{equation}
\E\bigl( \vert u_t(x)-u_t(y)\vert^{2k}\bigr) \le
C^{k^2}\biggl(\frac{|y-x|}{\varkappa}\biggr)^k \exp\biggl(\frac{B
k^3}{\varkappa}\biggr)
\end{equation}
for constants $B, C\in(0,\infty)$ that do not depend on $k$.
We apply an argument,
similar to one we used in the proof of Lemma~\ref{lemMC2},
in order to deduce that for simultaneously all intervals $I$ of length 1,
%
\begin{equation}\label{padiff}
\E\Bigl({\sup_{x,x' \in I\dvtx |x-x'|\le1}}|u_t(x)-u_t(x')|^{2k}\Bigr)
\le
\frac{C_1^{k^2}{\mathrm e}^{C_2k^3/\varkappa}}{\varkappa^k}
\end{equation}
for constants $C_1, C_2\in(0,\infty)$ that do not depend on $k$ or
$\varkappa$.
Now, we follow the proof of Theorem~\ref{thboundedsigma}
and partition $[-R,R]$ into intervals of length $1$. Let $b > 0$
to deduce the following:
%
\begin{equation}
\P\bigl\{ u_t^{*}(R) > 2 {\mathrm e}^{b(\log R)^{2/3}/\varkappa
^{1/3}}\bigr\}
\le\mathcal{T}_1+\mathcal{T}_2,
\end{equation}
where
%
\begin{equation}
\mathcal{T}_1 := \P\Bigl\{
\max_{1\le j \le2R} u_t(x_j) > {\mathrm e}^{b(\log R)^{2/3}/\varkappa
^{1/3}}\Bigr\}
\end{equation}
and
%
\begin{equation}\quad
\mathcal{T}_2 := \P\Bigl\{
{\max_{1\le j \le2R} \sup_{x\in(x_j,x_{j+1})}}|u_t(x)-u_t(x_j)| >
{\mathrm e}^{b(\log R)^{2/3}/\varkappa^{1/3}}\Bigr\}.
\end{equation}
[Compare with (\ref{eqprustar}).]

On one hand, Lemma~\ref{lempatails} implies that
%
\begin{equation}\label{pat1}
\mathcal{T}_1 \le2R \cdot
\P\bigl\{ u_t(x_j) >{\mathrm e}^{b(\log R)^{2/3}/\varkappa^{1/3}}
\bigr\}
\le\frac{2 c_3 R}{R^{c_4b^{3/2}}}
\end{equation}
for some constants $c_3,c_4>0$. On the other hand (\ref{padiff}) and
Chebyshev's inequality imply that
%
\begin{eqnarray}\label{pat2}
\mathcal{T}_2 &\le& 2R \P\Bigl\{
{\sup_{x,x'\in I\dvtx |x-x'|\le1}}|u_t(x)-u_t(x')|
\ge{\mathrm e}^{b(\log
R)^{2/3}/\varkappa^{1/3}}\Bigr\}\nonumber\\[-8pt]\\[-8pt]
&\le& \frac{2RC_1^{k^2}{\mathrm e}^{C_2k^3/\varkappa}}{\varkappa^k
{\mathrm e}^{2kb(\log R)^{2/3}/\varkappa^{1/3}}}.\nonumber
\end{eqnarray}
%
Now we choose $k:=\lceil\varkappa^{1/3} (\log R)^{1/3}
\rceil$
in order to obtain
$\mathcal{T}_2\le\mathrm{const}\cdot R^{1+C_2-2b}$ where the constant
depends on $\varkappa$. With these choices of parameters
we deduce from
(\ref{pat1}) and (\ref{pat2}) that if $b$ were sufficiently large, then
%
\begin{equation}
\sum_{R=1}^{\infty} \P\bigl\{
u_t^{*}(R) > 2 {\mathrm e}^{b(\log R)^{2/3}/\varkappa^{1/3}}\bigr\}
<\infty.
\end{equation}
The Borel--Cantelli lemma and a monotonicity argument
together complete the proof.
\end{pf*}

\begin{appendix}

\section{Walsh stochastic integrals}\label{appSI}
Throughout this Appendix, $(\Omega,\F,\P)$ denotes
(as is usual) the underlying probability space.
We state and prove some elementary properties of Walsh
stochastic integrals~\cite{Walsh}.

Let $\mathcal{L}^d$ denote the collection of all Borel-measurable sets
in $\R^d$ that have finite $d$-dimensional Lebesgue measure.
(We could work with Lebesgue-measurable sets, also.)

Let us follow Walsh~\cite{Walsh}
and define for every $t>0$ and $A\in\mathcal{L}^d$
the random field
%
\begin{equation}
W_t(A) := \int_{[0,t]\times A} W(\d s \,\d y).
\end{equation}
The preceding stochastic integral is defined in the same sense as Wiener.

Let $\F_t(A)$ denote the sigma-algebra generated by all random
variables of the form
%
\begin{equation}
\{ W_s(B)\dvtx s\in(0,t], B\in\mathcal{L}^d, B\subseteq
A\}.
\end{equation}
We may assume without loss of generality that, for all
$A\in\mathcal{L}^d$, $\{\F_t(A)\}_{t>0}$ is a right-continuous
$\P$-complete filtration (i.e., satisfies the ``usual hypotheses''
of Dellacherie and Meyer~\cite{DM}). Otherwise, we augment $\{\F
_t(A)\}_{t>0}$
in the usual way. Let
%
\begin{equation}
\F_t := \bigvee_{A\in\mathcal{L}^d} \F_t(A)\qquad
(t>0).
\end{equation}
Let $\mathcal{P}$ denote the collection of all
processes that are predictable with respect to $\{\F_t\}_{t>0}$.
The elements of $\mathcal{P}$ are precisely the ``predictable
random fields'' of Walsh~\cite{Walsh}.

For us, the elements of $\mathcal{P}$ are of interest because if
$Z\in\mathcal{P}$ and
%
\begin{equation}
\|Z\|^2_{L^2(\R_+\times\R^d\times\Omega)}
:=\E\int_0^\infty\d t\int_{\R^d}\d x\, [ Z_t(x)
]^2<\infty,
\end{equation}
then the Walsh stochastic integral $I_t:=\int_{[0,t]\times\R}
Z_s(y) W(\d s \,\d y)$ is defined properly, and has good mathematical
properties. Chief among those good properties are the following:
$\{I_t\}_{t>0}$ is a\vspace*{1pt} continuous mean-zero $L^2$-martingale with
quadratic variation $\langle I\rangle_t :=
\int_0^t\d s\int_{-\infty}^\infty\d y\, [Z_s(y)]^2$.

Let us define $\mathcal{P}(A)$ to be the collection of all processes
that are predictable with respect to $\{\F_t(A)\}_{t>0}$. Clearly,
$\mathcal{P}(A)\subseteq\mathcal{P}$ for all $A\in\mathcal{L}^d$.
%
\begin{proposition}\label{prmeasWalshint}
If\vspace*{1pt} $Z\in\mathcal{P}(A)$ for some $A\in\mathcal{L}^d$
and $\|Z\|_{L^2(\R_+\times\R^d\times\Omega)}<\infty$,
then the martingale defined by
$J_t:= \int_{[0,t]\times A} Z_s(y) W(\d s \,\d y)$
is in $\mathcal{P}(A)$.
\end{proposition}
\begin{pf}
It suffices to prove this for a random field $Z$
that has the form
%
\begin{equation}
Z_s(y) (\omega) = \one_{[a,b]}(s) X(\omega)\one_A(y)\qquad
(s>0, y\in\R, \omega\in\Omega)\vadjust{\goodbreak}
\end{equation}
where $0\le a<b$, and $X$ is a bounded $\F_a(A)$-measurable
random variable. But in that case,
$J_t(\omega)=X(\omega)\cdot\int_{[0,t]\cap[a,b]\times A}W(\d s \,\d y)$,
whence the result follows easily from the easy-to-check
fact that the stochastic process defined by
$I_t:=\int_{[0,t]\cap[a,b]\times A}W(\d s \,\d y)$
is continuous (up to a modification). The latter
assertion follows from
the Kolmogorov continuity theorem; namely, we check first that
$\E(|I_t-I_r|^2)=|A| \cdot|t-r|$, where $|A|$ denotes the
Lebesgue measure of $A$. Then use the fact, valid for all Gaussian
random variables including $I_t-I_r$, that
$\E(|I_t-I_r|^k)=\mathrm{const}\cdot\{ \E(|I_t-I_r|^2)\}^{k/2}$
for all $k\ge2$.
\end{pf}

Proposition~\ref{prmeasWalshint} is a small variation on
Walsh's original construction of his stochastic integrals.
We need this minor variation for the following reason:

\begin{corollary}\label{comeasWalshint}
Let $A^{(1)},\ldots,A^{(N)}$ be fixed and nonrandom
disjoint elements of
$\mathcal{L}^d$. If $Z^{(1)},\ldots,Z^{(N)}$ are, respectively, in
$\mathcal{P}(A^{(1)}),\ldots,\mathcal{P}(A^{(N)})$ and
$\|Z^{(j)}\|_{L^2(\R_+\times\R^d\times\Omega)}<\infty$
for all $j=1,\ldots,N$, then $J^{(1)},\ldots,J^{(N)}$
are independent processes, where
%
\begin{equation}
J^{(j)}_t := \int_{[0,t]\times A_j} Z_s(y)
W(\d s \,\d y)\qquad (j=1,\ldots,N, t>0).
\end{equation}
\end{corollary}
\begin{pf}
Owing to Proposition~\ref{prmeasWalshint}, it suffices to
prove that if some sequence of random fields $X^{(1)},\ldots,X^{(N)}$
satisfies $X^{(j)}\in\mathcal{P}(A^{(j)})$
($j=1,\ldots,N$), then $X^{(1)},\ldots,X^{(N)}$ are independent.
It suffices to prove this in the case that the $X^{(j)}$'s
are simple predictable processes; that is, in the case that
%
\begin{equation}
X^{(j)}_t(\omega) = \one_{[a_j,b_j]}(s) Y_j(\omega) \one_{A^{(j)}}(y),
\end{equation}
where $0<a_j<b_j$ and $Y_j$ is a bounded $\F_{a_j}(A^{(j)})$-measurable
random variable. In turn, we may restrict attention to $Y_j$'s that
have the form
%
\begin{equation}
Y_j(\omega) := \varphi_j\biggl(\int_{[\alpha_j,\beta_j]
\times A^{(j)}} W(\d s \,\d y)\biggr),
\end{equation}
where $0<\alpha_j<\beta_j \le a_j$ and $\varphi_j\dvtx\R\to\R$ is
bounded and Borel measurable. But the assertion is now clear, since
$Y_1,\ldots,Y_N$ are manifestly independent. In order to see this
we need only verify that the covariance between
$\int_{[\alpha_j,\beta_j]\times A^{(j)}}W(\d s \,\d y)$
and $\int_{[\alpha_k,\beta_k]\times A^{(k)}}W(\d s \,\d y)$
is zero when $j\neq k$; and this is a ready consequence of
the fact that $A^{(j)}\cap A^{(k)}=\varnothing$ when
$j\neq k$.
\end{pf}

\section{Some final remarks}\label{secasinterm}
Recall that $\|U\|_k$ denotes the usual $L^k(\P)$-norm
of a random variable $U$ for all $k\in(0,\infty)$.
According to Lemma~\ref{lemmomentsU},
%
\begin{equation}
\limsup_{t\to\infty} \frac1t \log\sup_{x\in\R}
\E(|u_t(x)|^k)
\le aCk^3 \qquad\mbox{if $k\ge2$}.
\end{equation}
This and Jensen's inequality together imply that
%
\begin{equation}
\overline\gamma(\nu) := \limsup_{t\to\infty}
\frac1t \log\sup_{x\in\R}\E(|u_t(x)|^\nu
)<\infty\qquad
\mbox{for all $\nu>0$}.
\end{equation}
And Chebyshev's inequality implies that for all $\nu>0$,
%
\begin{eqnarray}
\sup_{x\in\R}\P\{ u_t(x)\ge{\mathrm e}^{-qt}\} &\le&
\exp\biggl(\nu t\biggl[q + \frac{1}{\nu t}
\log\sup_{x\in\R}\E(|
u_t(x)|^\nu)\biggr]\biggr)\nonumber\\[-8pt]\\[-8pt]
&=& \exp\biggl( \nu t \biggl[ q +
\frac{\overline{\gamma}(\nu)}{\nu}+o(1)\biggr]\biggr)\qquad
(t\to\infty).\nonumber
\end{eqnarray}
Because $u_t(x)\ge0$, it follows that $\nu\mapsto\overline{\gamma
}(\nu
)/\nu$
is nondecreasing on $(0,\infty)$, whence
%
\begin{equation}\label{ell}
\ell:= \lim_{\nu\downarrow0}\frac{\overline{\gamma}(\nu)}{\nu}
=\inf_{\nu>0}\frac{\overline{\gamma}(\nu)}{\nu}
\mbox{ exists and is finite}.
\end{equation}
Therefore, in particular,
%
\begin{equation}
\limsup_{t\to\infty}
\frac1t \log\sup_{x\in\R}\P\{ u_t(x)\ge{\mathrm e}^{-qt}
\}<0\qquad\mbox{for all $q\in(-\infty,-\ell)$}.
\end{equation}
Now consider the case that $\sigma(0)=0$, and recall that
in that case
$u_0(x)\ge0$ for all $x\in\R$. Mueller's comparison principle
tells us that
$u_t(x)\ge0$ a.s. for all $t\ge0$ and $x\in\R$, whence
it follows that $\|u_t(x)\|_1=\E[ u_t(x)]=(p_t*u_0)(x)$ is bounded
in $t$. This shows that $\overline{\gamma}(1) =0$, and hence
$\ell\le0$. We have proved the following:
%
\begin{proposition}\label{prBG}
If $\sigma(0)=0$, then there exists $q\ge0$ such that
%
\begin{equation}\label{preqBG}
\frac1t \log u_t(x) \le-q + o_{_{\P}}(1) \qquad\mbox{as $t\to
\infty$,
for every $x\in\R$},
\end{equation}
where $o_{_{\P}}(1)$ is a term that converges to zero in probability as
$t\to\infty$.
\end{proposition}

Bertini and Giacomin~\cite{BG} have studied the case that
$\sigma(x)=cx$ and have shown that, in that case,
there exists a special choice of $u_0$ such that for all
compactly supported probability densities $\psi\in C^\infty(\R)$,
%
\begin{equation}\label{eqBG}
\frac1t \log\int_{-\infty}^\infty u_t(x)
\psi(x) \,\d x = -\frac{c^2}{24\varkappa} + o_{_{\P}}(1)
\qquad\mbox{as $t\to\infty$}.
\end{equation}
Equation (\ref{eqBG}) and more generally Proposition
\ref{prBG} show that the typical behavior of the sample function
of the solution to (\ref{eqmain}) is subexponential in time,
as one might expect from the unforced linear heat equation.
And yet, it frequently is the case that $u_t(x)$ grows in time exponentially
rapidly in $L^k(\P)$ for $k\ge2$ \mbox{\cite{BC,CM,FK}}. This phenomenon
is further evidence of physical intermittency in the sort of systems
that are modeled by (\ref{eqmain}).

Standard predictions suggest that the typical behavior of $u_t(x)$
(in this and related models) is that it decays exponentially rapidly
with time.\vadjust{\goodbreak} [Equation (\ref{eqBG}) is proof of this fact in one special
case.] In other words, one might expect that typically $q>0$.
We are not able to resolve this matter here, and therefore ask
the following questions:\vspace*{8pt}

\textit{Open problem} 1. Is $q> 0$ in Proposition~\ref{prBG}?
Equivalently, is $\ell<0$ in~(\ref{ell})?\vspace*{8pt}

\textit{Open problem} 2. Can the $o_{_{\P}}(1)$ in (\ref{preqBG})
be replaced by a term that converges almost surely to zero as $t\to
\infty$?
\end{appendix}



\printaddresses

\end{document}